\documentclass[11pt]{amsart}
\usepackage{times}      % use Times fonts for text
\usepackage{latexsym}   % use all LaTeX fonts
\usepackage{amssymb}    % use all AMS fonts
\usepackage{amsmath}    % include some AMS-LaTeX functionality
\usepackage{amsthm}     % AMS style theorem and proof environments
\usepackage{mathrsfs}   % special calligraphic characters like \mathscr{S}
\usepackage{euscript}   % euscript package
\usepackage{epsfig,graphics,color,subfigure}             
\usepackage[all]{xy}    % Xy-Pic for diagrams
\newcommand{\ra}{\rightarrow}
\newcommand{\R}{\mathbb{R}}

\newcommand{\C}{\mathbb{C}}

\renewcommand{\epsilon}{\varepsilon}

\renewcommand{\qed}{\square}
\renewcommand{\"}[1]{\mathaccent 127 #1}
\renewcommand{\phi}{\varphi}

\newcommand{\cc}{d_{cc}}

\renewcommand{\"}[1]{\mathaccent 127 #1}
\renewcommand{\phi}{\varphi}

\newcommand{\ip}{<\cdot,\cdot>}
\newcommand{\pf}{\noindent {\em Proof: }}
 
\newtheorem{Pro}{Proposition}[section]
\newtheorem{Lem}[Pro]{Lemma}

\newtheorem{Thm}[Pro]{Theorem}
\newtheorem{Qu}{Question}
\newtheorem{MThm}{Theorem}

\newtheorem{Rem}{Remark}
\newtheorem{Def}[Pro]{Definition}

\newcommand{\hook}{\lfloor}
\begin{document}
\newcommand{\wra}{\rightharpoonup}
\renewcommand{\theenumi}{\roman{enumi}}
\newtheorem{Prob}[Pro]{Problem}

\title{Minimal surfaces in the Heisenberg group}
\author{Scott D. Pauls}
\thanks{The author is partially support by NSF
    grant DMS-9971563}
\address{Dartmouth College, Hanover, NH, 03755}
\email{scott.pauls@dartmouth.edu}
\keywords{Minimal surfaces, Carnot groups, Heisenberg group, 1-Laplacian}

\begin{abstract}

We investigate the minimal surface problem in the three dimensional Heisenberg
group, $H$, equipped with its standard Carnot-Carath\'eodory metric.
Using a particular surface measure, we characterize minimal surfaces
in terms of a sub-elliptic partial differential equation and prove an
existence result for the Plateau problem in this setting.  Further,
we provide a link between our minimal surfaces and Riemannian constant
mean curvature surfaces in $H$ equipped with different Riemannian metrics
approximating the Carnot-Carath\'eodory metric.  We generate a large
library of examples of minimal surfaces and use these to show that the
solution to the Dirichlet problem need not be unique.  Moreover, we show
that the minimal surfaces we construct are in fact X-minimal surfaces
in the sense of Garofalo and Nhieu.       
\end{abstract}

\maketitle
\section{Introduction}
The examination of minimal surfaces in various settings has a long and
rich history.   An exploration of minimal surfaces in
$\R^3$, $\R^n$, and Riemannian manifolds led to beautiful and amazing
connections to complex analysis, harmonic mappings, and many other
diverse settings.  These connections led to the solution of many
problems in minimal surface theory, providing many different methods
for the construction of examples and different proofs of existence and
uniqueness results.  For example,  variational techniques were used to
examine minimal surfaces in terms of a partial differential equation
and to address existence and uniqueness of solutions with prescribed
boundary data via an examination of the partial differential equation.  These
techniques led to very general existence and uniqueness results for
minimal surfaces in many different settings.  However, even with these decades of investigation and
exploration, the so-called {\em Plateau problem} remains compelling and not yet
completely understood in many settings: 

\vspace{.2in}

\noindent
{\bf Plateau Problem: }{\em For a given curve, $\gamma$, can one find a surface of
least area spanning $\gamma$? }

\vspace{.2in}

The purpose of this paper is to continue an investigation of this
basic question in the setting of the three dimensional Heisenberg
group, $H$, equipped with a Carnot-Carath\'eodory (CC) metric. To
define such a metric, we first must describe some of the features of $H$.
$H$ is a Lie group defined by a Lie algebra $\mathfrak{h}$ generated
by three vector fields, $\{X,Y,Z\}$ with one nontrivial bracket
relation, namely $[X,Y]=Z$.  Moreover, we use this to describe the
grading of $\mathfrak{h}$, namely that \[ \mathfrak{h} = \mathcal{V}
\oplus \mathcal{V}_2\] where $\mathcal{V} = span\{X,Y\}$ and
$\mathcal{V}_2 = span\{ Z\}$.  While there are many presentations of
the Heisenberg group, for the purposes of this paper, we use an
identification of $H$ with $\R^3$ where the vector fields $\{X,Y,Z\}$
are given by 
\begin{equation*}
\begin{split}
X&= \partial_x -y \partial_z\\
Y&= \partial_y+x \partial_z\\
Z&= 2\partial_z
\end{split}
\end{equation*}
where $\{\partial_x,\partial_y,\partial_z\}$ are the standard basis
vector fields in $\R^3$.  Using left translation by group
elements, we can think of $\mathcal{V}$ and $\mathcal{V}_2$ as left
invariant subbundles of the tangent bundle of $H$.  Abusing notation,
we will refer to both the vector space and subbundle by the same
symbol.  Throughout this paper, we often call 
$\mathcal{V}$ the {\em bottom of the grading} of $\mathfrak{h}$ and,
if a path $p \subset H$ has tangent vector in $\mathcal{V}$ at every
point, we call $p$ a {\em horizontal curve}.  Now, fixing a left
invariant inner product, $\ip$, on $\mathcal{V}$ making $\{X,Y\}$ an
orthonormal basis, we define
the CC metric as a path metric on $H$ as follows:

\[ \cc(g,h) = \inf \left \{ \int_I <p'(t),p'(t)>^\frac{1}{2} \; dt
\bigg | p(0)=g, \;
p(1)=h, \; \text{ and $p$ is horizontal} \right \}\]

It follows from the fact that $\mathcal{V}$ is bracket generating that
$\cc(g,h) < \infty$ for any $g,h \in H$.  Notice that, by definition,
the CC metric is left invariant.  Further the CC
metric is fractal in the sense that, while the topological dimension
of $H$ is three, the Hausdorff dimension of $H$, calculated with
respect to the CC distance function, is four.  In this setting, we investigate
a basic version of the question above:

\vspace{.2in}

\noindent
{\bf Plateau Problem in $H$: }{\em Given a closed curve $\gamma
\in H$, can we find a topologically two dimensional surface $S \subset
H$ spanning $\gamma$ which minimizes an appropriate surface measure? }

\vspace{.2in}

The singularity of the CC metric, as illustrated by the disconnect between the topological and Hausdorff
dimension mentioned above, makes the notion of ``appropriate surface
measure'' somewhat hard to specify - there are several ``natural''
candidates for such a measure.  For example, in \cite{Gromov:CC}, among many other
investigations, Gromov addresses such so-called filling problems,
although focusing more on isoperimetric type inequalities and
questions (see, for example, section 0.7B).  He shows that one natural measure for the generic two
dimensional surface in $H$ is the three dimensional Hausdorff measure
associated to $\cc$ as every $C^1$ surface which has topological
dimension two has Hausdorff dimension three (this is an application of
results in section 0.6A.  Also, the reader should see chapter 2 section 1 for
a more general discussion).  In a different line of investigation, using a measure theoretic
perimeter measure, $\mathscr{P}$, Garofalo and Nhieu (\cite{GarNh}) gave a beautiful
and delicate
existence result for minimal surfaces in a much more general setting than the Heisenberg
group.  For the sake of simplicity, we state their theorem only in the
case of $H$ and precise definitions for the necessary objects appear
in section \ref{back}.  To differentiate between their result and
others, we follow their convention and 
call their minimal surfaces {\em X-minimal surfaces}.  

\begin{Thm}\label{GN}(\cite{GarNh})   Given a bounded open set $O
\subset H$ and an X-Caccioppoli set $L \subset H$, there exists an
X-minimal surface in $O$, i.e. an X-Caccioppoli set $E \subset H$
such that $\mathscr{P}(E) \le \mathscr{P}(F)$ for every set $F$ which coincides with
$L$ outside of $\Omega$.
\end{Thm}

To interpret this theorem in the language of the question posed above,
we must indicate where the curve $\gamma$ and the two dimensional
surface appear.  Garofalo and Nhieu's perimeter measure is an analogue
of DeGiorgi's perimeter measure in Euclidean space and it is a type of
area measure on the measure theoretic boundaries of open sets.  Thus,
to recover a ``curve'', we simply intersect the measure theoretic
boundaries of $L$ and $O$ in the theorem above.  The surface
spanning this curve is then the boundary of $E$ inside $O$.
Notice that Garofalo and Nhieu's work addresses only existence and not
regularity hence the measure theoretic nature of the result.  

We begin our investigation from an a priori different direction.
We will consider surfaces bounding a specific
curve which minimize an energy based on the three dimensional spherical Hausdorff measure.  We use
standard variational techniques to find a partial differential
equation characterizing these minimal surfaces in the Heisenberg group.
If the surface, $S$, is given by
$F(x,y,z)=f(x,y)-z=0$ for $(x,y) \in \Omega$, a region in the
$xy$-plane (here we are using the indentification of $H$ and $\R^3$,
then it is known (see \cite{Pansu} or \cite{CCcalc}) that, up to a
normalization, the
spherical Hausdorff measure takes the form: 
\[\mathscr{H}^3_{cc}(S) = \int_\Omega |\nabla_0 F| dA \]
where $\nabla_0$ is the so-called {\em horizontal gradient operator}
given by 
\[\nabla_0 F = (XF,YF)\]
Using this as the basis for our variational setup, we define our
energy function to be \[E(\cdot)= \int_\Omega |\nabla_0 \cdot| dA\]
Under these assumption, 
 we can
characterize nonparametric minimal surfaces as follows:  if the minimal surface is
given as the level set $F(x,y,z)=f(x,y)-z=0$, then it satisfies the
partial differential equation:
\begin{equation}\label{inteq}
 \nabla_0 \cdot \frac{\nabla_0 F}{|\nabla_0 F|}=0
\end{equation}

This ``minimal surface equation'' is a subelliptic partial differential equation and one can not
immediately conclude existence or uniqueness results from an
examination of the defining energy functional.  To investigate the standard
questions of existence and uniqueness, and to relate these minimal
surfaces to the X-minimal surfaces of Garofalo and Nhieu, we rely on
approximation of $(H,\cc)$ by Riemannian manifolds. 

It is well know (see, for example \cite{GLP}, \cite{Gromov:polygrowth},
and \cite{Pansu}) that $(H,\cc)$ can be realized as a limit of
dilated Riemannian manifolds.  Specifically, we first fix a left
invariant Riemannian metric, $g_1$, on $H$ matching $\ip$ on
$\mathcal{V}$ and making $\{X,Y,Z\}$ an orthonormal basis for $\mathfrak{h}$.  Then, we define a dilation
map, $h_\lambda:H \ra H$ by 
\[ h_\lambda(e^{aX+bY+cZ}) = e^{a\lambda X+b\lambda Y+c\lambda^2Z}\]
we define a new Riemannian metric on $H$ for every $\lambda$ by
\[ g_\lambda = \frac{1}{\lambda} h_\lambda^* g_1 \]
Then, if we denote the distance function associated to $g_\lambda$ by
$d_\lambda$, the sequence of metric spaces $(H,g_\lambda)$ converges
to $(H,\cc)$ in the Gromov-Hausdorff topology as $\lambda \ra
\infty$.  

Using this convergence, we attack the minimal surface problem in
$(H,\cc)$ by considering sequences of surfaces, $\{S_\lambda\}$ where 
$S_\lambda \subset (H,g_\lambda)$.  In particular, we can use sequences of Riemannian
minimal surfaces and sequences of Riemannian constant curvature
surfaces to construct solutions to (\ref{inteq}), i.e. minimizers of our
energy functional.  This allows us to guarantee existence, at least
weakly, for a special class of curves:

\begin{MThm}\label{I1}
 Let $\Gamma$ be a closed curve in $H$ satisfying the bounded
  slope condition and which is the graph
  of a function $\phi\in C^{2,\alpha}(\R)$ over a curve $\gamma \in C^{2,\alpha}(\R^2)$ which bounds a
  region $\Omega$.  Then, there exists $u \in W^{1,p}(\Omega) \cap C^0(\Omega)$ so that
  $u|_\gamma = \phi$, $u$ is a weak solution to (\ref{inteq}) on
  $\Omega$ and $u$ minimizes $E(\cdot)$ on $\Omega$.  Moreover, there
  exists a sequence of functions with the same regulatity as $\Gamma$,
  $\{u_{\lambda_n}\}$, such that $H_{\lambda_n}(u_{\lambda_n})=0$ and
  $u_{\lambda_n} \ra u$ in $W^{1,p}(\Omega)$.
\end{MThm}

This theorem is proved by showing that the sequence of minimal
surfaces with boundary $\Gamma \subset (H,g_\lambda)$ have a
convergent subsequence as $\lambda \ra \infty$ and then showing that
the limit is a graph satisfying (\ref{inteq}).  Basically, this allows
us to say that, upon fixing boundary conditions, Riemannian minimal surfaces
subconverge to minimal surfaces in the Carnot setting.  We note that
the bounded slope condition is not optimal - any condition on the
curve $\Gamma$ which ensures existence of minimal surfaces in
$(H,g_\lambda)$ is sufficient.  While this
allows us to conclude an existence result, it says nothing about
uniqueness.  In fact, the next theorem gives strong evidence for
nonuniqueness, by showing that limits of families of constant
mean curvature surfaces (if such a limit exists) are minimal surfaces
in $(H,\cc)$ as well.

\begin{MThm}\label{I2}
  Let $\Gamma$ be a closed curve with the
bounded slope condition in $H$ which is the graph
  of a function $\phi \in C^{2,\alpha}(\R)$ over a curve $\gamma \in C^{2,\alpha}(\R^2)$ which bounds a
  region $\Omega$.  Suppose $u_n$ are graphs in $(H,d_{\lambda_n})$
  spanning $\Gamma$ which have mean curvature given by the function
  $\kappa_n$.  If $\kappa_n \ra 0$ as $n \ra \infty$
  and $u_n$ converge to a graph $u$ then $u$ satisfies (\ref{inteq}) and
  is a minimizer of $E(\cdot)$.  
\end{MThm}

As with the previous theorem, this theorem requires the bounded slope condition to ensure existence
of minimal surfaces with the same boundary data in $(H,g_\lambda)$.
With regard to the uniqueness question, one now suspects that there
may be many different solutions to the Dirichlet problem in this setting
generated by taking different convergent subsequences.  This suspicion
is confirmed in section \ref{unique} where we construct an explicit
example where there are at least two minimal surfaces spanning the
same curve.  Hence we conclude:

\begin{MThm}\label{I3}  Let $\Gamma$ be a closed $C^2$ curve in $H$ which is the
  graph of a function $\phi$ over a curve $\gamma \in \R^2$ which
  bounds a region $\Omega$.  Then the solution to the Dirichlet problem
  for $\Gamma$ need not be unique.
\end{MThm}

This theorem is a consequence of the work in section \ref{egs}, where
we construct a huge number of solutions to equation $\ref{inteq}$.  We
note that, in addition to being examples of the minimal surfaces
described above, these are solutions to the $1$-subLaplacian equation
\ref{inteq} on $H$ as well as solutions to the least horizontal gradient problem on
$\R^3$.  Next, we give a compilation of all of the examples, shown via
various techniques, in section \ref{egs}.

\begin{MThm}\label{Ieg}
Given $a,b,c \in \R$, $\alpha \in [-1,1]$
  and $g \in W^{1,p}(\R)$, the following graphs satisfy (\ref{inteq})
  weakly:
\begin{enumerate}
\item \[z=a x+b y+c\]
\item \[z=a \theta\] (using cylindrical coordinates on $\R^3$)
\item \[z= \pm \left ( \frac{\sqrt{b r^2-1}}{b} - a
      \tan^{-1} \left ( \frac{1}{\sqrt{b r^2-1}} \right ) \right )
    + c + a \theta\]  Again, these use cylindrical
    coordinates.
\item \[z=xy+g(y)\]
\item \[ z = \frac{\alpha}{\sqrt{1-\alpha^2}} x^2+xy+g\left(
    y+\frac{\alpha}{\sqrt{1-\alpha^2}}x\right)\]
\item \[ z = -\frac{\sqrt{1-\alpha^2}}{\alpha} x^2+xy+g\left(
    y-\frac{\sqrt{1-\alpha^2}}{\alpha}x\right)\]
\item \[z=g \left ( \frac{y}{x}\right)\]
\item \[z=\frac{a}{b+x}x^2+\frac{b^2}{b+x}xy+g \left
    ( \frac{y-a}{b+x} \right )\]
\end{enumerate}
\end{MThm}
Although these results are promising,we began with a different measure than in Garofalo and
Nhieu's on X-minimal surfaces and, a priori, our minimal surfaces have
no relation to X-minimal surfaces.  However, as mentioned above, the
Riemannian approximation spaces allow us to demonstrate a link between
the two:

\begin{MThm}\label{I4}
The minimal surfaces described in theorems \ref{I1} and \ref{I2} are
X-minimal as well.
\end{MThm}

The proof of this theorem rests heavily on standard elliptic theory in
Riemannian 
manifolds coupled with Garofalo and Nhieu's results and results on
X-Caccioppoli sets due to Franchi, Serapioni and Serra Cassano in
\cite{FSSC} and \cite{FSSC2}.  In particular, the results in
\cite{FSSC2} include a form of the implicit function theorem which
allows us to compare the graphs of minimal surfaces with the graphs of
X-Caccioppoli sets.  

After the completion of this work, the author learned through a
personal communication from N. Garofalo that he, D. Danielli, and
D.M. Nhieu (\cite{DGN}) had explored the theory of X-minimal surfaces
further and had independently arrived at similar results using
somewhat different techniques.  In particular, they find the same
characterization of X-minimal surfaces given by equation \ref{inteq}
and also address aspects of the Bernstein problem and the construction
of examples.  In addition, they provide numerous results in other
directions.  

To conclude the introduction, we summarize the contents of each section.  Section \ref{back}
reviews many of the definitions above as well as other necessary
definitions and results from various sources.  Section \ref{char}
provides the variational analysis of the energy functional described
above and proves theorems \ref{I1}, \ref{I2} and \ref{I4}.  In section
\ref{egs}, we construct, through a variety of methods, the examples
of minimal surfaces in $(H,\cc)$ stated in theorem \ref{Ieg}.  Section \ref{cons} points out
several consequences of the results in section \ref{char} and the
examples in section \ref{egs}.  In particular, it describes the
example proving theorem \ref{I3} and discusses the implications of the
examples with reference to the classical Bernstein problem.  The end
of this section discusses some conclusions and avenues of continued
exploration.

We wish to thank the referee for many helpful comments and suggestions.
\section{Further definitions and a recounting of known
  results}\label{back}

For the convenience of the reader, we begin this section by collecting the definitions embedded in the
introduction:
\begin{itemize}
\item $H$ denotes the three dimensional Heisenberg group.  It is
  associated to the Lie algebra $\mathfrak{h}= span\{X,Y,Z\}$ with one
  nontrivial bracket, $[X,Y]=Z$.
\item Using the exponential map, we will often identify $H$ with $\R^3$
  using $\{X,Y,Z\}$ as the standard coordinates.  We will often
  describe graphs in $H$ and will use the convention that the graph of
  a function
  $g:\R^2 \ra \R$ is given by the set $(x,y,g(x,y))$ unless otherwise
  specified.  
\item $\mathfrak{h}$ has a grading given by
  \[\mathfrak{h} = \mathcal{V} \oplus \mathcal{V}_2\]
where $\mathcal{V} = span\{X,Y\}$ and $\mathcal{V}_2=span\{Z\}$.
\item $\mathcal{V}$ is thought of both as a subspace of $\mathfrak{h}$
  and as a left invariant vector bundle on $H$ via left translation.
  We will refer to $\mathcal{V}$ as the {\em bottom of the grading} or
  as the {\em horizontal bundle} on $H$.
\item Paths whose tangent vectors are always horizontal are called
  {\em horizontal paths}.
\item $\ip$ denotes the standard inner product on $\mathcal{V}$,
  i.e. the inner product which makes $\{X,Y\}$ an orthonormal basis.
\item $\nabla_0$ denotes the horizontal gradient operator and is
  defined by \[\nabla_0f=(Xf,Yf)\]
\item Given a $C^1$ surface, $S$, in $H$, we denote by $N_0$ the {\em
    horizontal normal vector} given by projecting the usual normal
  vector to the horizontal bundle at each point.  If we wish to
  emphasize the dependence on $S$, we will write $N_0(S)$.  Similarly,
  to emphasize the dependence of the normal on the base point, we will
  write $N_0(p)$ where $p$ is a point on the surface.  
\item A point $p$ on a surface $S$ is called a {\em characteristic
    point} if $N_0$ vanishes at $p$.
\item Let $n_0= \frac{N_0}{|N_0|}$ be the unit horizontal normal.
  Given a surface $S$, we often think of $n_0$ as a map from $S$ to
  the circle.  In this case, we call $n_0$ the {\em horizontal Gauss map}. 
\item The standard Carnot-Carath\'eodory metric on $H$ is given by
\[ \cc(g,h) = \inf \left \{ \int_I <p'(t),p'(t)>^\frac{1}{2} \; dt
\bigg | p(0)=g, \;
p(1)=h, \; \text{ and $p$ is horizontal} \right \}\]
for $g,h \in H$.
\item $\mathscr{H}^k_{cc}$ denotes the $k$-dimensional spherical
  Hausdorff measure constructed with respect to the
  Carnot-Carath\'eodory metric.
\item We denote the open ball of radius $r$ around a point $p \in H$ with respect to
  $\cc$ by $B_{cc}(p,r)$.
\item $g_1$ is the left invariant Riemannian metric on $H$ which makes
  $\{X,Y,Z\}$ an orthonormal basis at each point.
\item $h_\lambda:H \ra H$ is a dilation map defined by
\[ h_\lambda(e^{aX+bY+cZ}) = e^{a\lambda X+b\lambda Y+c\lambda^2Z}\] 
\item $g_\lambda$ is a dilated Riemannian metric on $H$ given by
 \[g_\lambda =\frac{1}{\lambda}h_\lambda^*g_1\]
Roughly, $g_\lambda$ measures $X$ and $Y$ directions the same way as
$g_1$ but changes the length of $Z$ by a factor of $\lambda$.
\item For a given surface, $N_\lambda$ denotes the normal to the
  surface computed with respect to $g_\lambda$.
\item If $\Omega \subset \R^2$, we denote by $W^{1,p}(\Omega)$ and
  $W^{1,p}_0(\Omega)$ the
  usual Sobolev spaces.
\item Throughout the paper, we will use the ``big O'' and ``little o''
  notation to describe the rates of decay of various functions.
  However, we will use the notation to emphasize the dependence on
  certain variables.  For example, $h(\epsilon,x,y,z)=o_\epsilon(1)$ means
  that \[\lim_{\epsilon \ra 0} h(\epsilon,x,y,z) =0 \]
\end{itemize}

\subsection{Surface measures}

Next, we review in more detail many of the constructions and theorems
mentioned in the introduction.  We begin with a description of
$\mathscr{H}^3_{cc}$ for smooth surfaces in $H$.  As we will be
focusing on graphs over portions of the $xy$-plane, we will state the
theorems in this setting.  The next theorem was first shown for $H$ in
\cite{Pansu} and was later extended to all Carnot groups by
J. Heinonen in \cite{CCcalc}.  Again, for the purposes of this paper,
we will state it only in $H$.  

\begin{Pro}\label{hausformula} (\cite{Pansu},\cite{CCcalc}) Let $S$ be
  a smooth surface. Then,
\begin{enumerate}
\item \[\mathscr{H}^3_{cc}(S) = \int_S \frac{|N_0|}{|N_1|} \; dA\] where $N_0$ is
  the horizontal normal described above, $N_1$ is the normal to $S$
  computed with respect to $g_1$ and $dA$ is the Riemannian
  area element induced by $g_1$.
\item In particular, if $S$ be a smooth
 surface given as the graph of $f:\Omega\subset \R^2 \ra \R$ in $H$.  Then, 
\[\mathscr{H}^3_{cc}(S) = \int_\Omega |\nabla_0 (f(x,y)-z)| \; dx dy\]
\end{enumerate}
\end{Pro}

In particular, this formula shows that $\mathscr{H}^3_{cc}$ is a
natural measure to use on topologically two dimensional surfaces in
$H$ as it is locally finite.  We will use this measure as the starting
point for our variational analysis.  One reason for this choice is
that, as demonstrated above, this measure has an extremely nice
presentation when the surface is smooth.  Moreover, as we will see
next, for smooth surfaces, this coincides with the perimeter measure
used by Garofalo and Nhieu.  To introduce the perimeter measure, we
first make several definitions following the notation of \cite{GarNh}
very closely.  Consider an open set $O \subset H$ where we think of
$H$ as identified with $\R^3$.  First, we define the weak horizontal Sobolev space of first order.  For $1 \le p
< \infty$, 

\[ \mathscr{L}^{1,p}(O) = \{f \in L^p(O) | Xf,Yf \in L^p(O)\}\]

In this case, if $f$ is not smooth, we understand $Xf$ and $Yf$ to be
 distributions.   
To continue, we next define 
\[ \mathscr{F}(O) = \{ \phi=(\phi_1,\phi_2) \in C^1_0(O;\R^2) | ||\phi||_\infty \le 1
\}\]
 
This allows us to define the {\bf $X$-variation} of a function $u \in
L^1_{loc}(O)$ by
\[Var_X(u;O)= \sup_{\phi \in \mathscr{F}(O)} \int_Ou(h)
  (X^*\phi_1(h)+Y^*\phi_2(h)) \;dV(h) \] where $X^*$ and $Y^*$ are the formal adjoints of $X$ and $Y$ and $dV$
is Haar measure on $H$.  Notice
that if $u$ is smooth and in $\mathscr{L}^{1,1}(O)$ then an
application of the divergence theorem yields that 
\[Var_X(u) = \int_O |\nabla_0(u)| \; dV \]

Finally, we can define the {\bf $X$-perimeter} of an open set $E \in H$
relative to the open set $O$ by

\[ \mathscr{P}(E;O) = Var_X(\chi_E,O)\] where $\chi_E$ is the characteristic function of $E$.  Naturally, we
say that a set has finite perimeter (with respect to $O$) if
$\mathscr{P}(E;O)<\infty$.  Moreover, a set $E$ is called an {\em
  $X$-Caccioppoli set} if $\mathscr{P}(E;O) < \infty$ for all open
sets $O \subset H$.   

The measure $\mathscr{P}$ was first introduced in \cite{CDG} and is the direct analogue of DeGiorgi's perimeter measure
for $\R^n$ introduced in \cite{DeGiorgi}.  Using this
notation, we can define an $X$-minimal surface.  

\begin{Def} Let $L\subset H$ be an $X$-Caccioppoli set and $O$ be an open set
  in $H$.  Then, a set $M$ is called an {\bf $X$-minimal surface} with
  respect to $O$ and $L$ if, for every set $S$ which coincides with
  $L$ outside of $O$, we have
\[\mathscr{P}(M;O) \le \mathscr{P}(S;O) \]
\end{Def}

One key component of the proof of theorem 1.1 is the lower
semicontinuity of the perimeter measure.

\begin{Lem}\label{semicont}  Fix and open set $O$ in $H$ and let $\{L_k\}$ be a sequence of
  $X$-Caccioppoli sets in $H$ converging to the set $L$ in the sense that $\chi(L_k) \ra \chi(L)$ in
  $L^1_{loc}$.  Then,
\[ \mathscr{P}(L;O)  \le \liminf_{k\ra \infty} \mathscr{P}(L_k;O)\]
\end{Lem}

A useful technical tool is to approximate sets of finite perimeter by
  smooth sets.  This can be done by adapting the standard
  mollification process.  We state this result (a version of results
  in \cite{FSSC96} and \cite{GarNh}) for the convenience of the reader:

\begin{Pro}\label{smoothapprox}  If $E$ is a bounded subset of $H$ with finite
  $X$-perimeter, then there exists a sequence, ${E_n}$, of smooth sets
  so that 
\[\lim_{n\ra \infty} \mathscr{P}(E_n) = \mathscr{P}(E) \]
\end{Pro}
 
  For the
purposes of this paper, we will be interested in surfaces given as
graphs over sets in the $xy$-plane (i.e. nonparametric minimal
surfaces).  To rectify this with Garofalo and Nhieu's notation, we
make the following conventions.  Given an open subset $\Omega
\subset \R^2$ and a curve $\Gamma \subset H$ given as the graph of a
function $\phi:\partial \Omega \ra \R$, we define $O$ to be an open
cylinder over $\Omega$, \[ O= \{ (x,y,z) \in H | (x,y) \in \Omega \}
\]
Further, we will define $L$ in the definition above by fixing a
function $f:\R^2 \ra \R$ whose graph spans $\Gamma$ and let 
\[ L = \{(x,y,z) \in H | z < f(x,y)\} \]
Similarly, when considering candidates for $X$-minimal surfaces using
this setup, we will specify the graph of a function $u:\Omega \ra \R$,
extend it to match $f$ outside $\Omega$ and define the open set $M$ by
\[ M = \{ (x,y,z) | z < u(x,y)\}\]

To summarize, we make the following definition:
\begin{Def}  When considering nonparametric $X$-minimal surfaces, given a
  set $\Omega \subset \R^2$, $\phi: \partial \Omega \ra \R$, $f:\R^2
  \ra \R$ and $u: \Omega \ra \R$ as above.  We define the perimeter of
  $u$ as
\[ \mathscr{P}(u) = \mathscr{P}(M;O) \]
and say that the graph of $u$ defines an $X$-minimal surface spanning
$\Gamma$ if $M$ is
an $X$-minimal surface with respect to $O$ and $L$.
\end{Def}

Similarly, we will use a similar convention when denoting the Hausdorff measure of
the graph of $u$ over $\Omega$,
\[\mathscr{H}^3_{cc}(u) \equiv \mathscr{H}^3_{cc}(u(\Omega))\]

Next, we state some results that make working with the perimeter measure more
tractable.  These results are due either to B. Franchi, R. Serapioni and F. Serra
Cassano (\cite{FSSC2}, see proposition 2.14) or L. Capogna,
D. Daneilli and N. Garofalo (\cite{CDG2} p. 211).  The reader should also
consult Z. Balough's paper \cite{Balogh} concerning various aspects of
the perimeter measure for surfaces as well as a study of different
surfaces measure in CC spaces by R. Monti and F. Serra Cassano (\cite{MSC}).

\begin{Pro}(\cite{CDG2},\cite{FSSC2})\label{measform}  Let $S$ be a $C^1$ surface
  bounding an open set $O$,
  then
\begin{enumerate}
\item \[\mathscr{P}(O) = \int_S \frac{|N_0|}{|N_E|} \; d\mathscr{H}^2_E\] where $N_0$ is the
  horizontal normal, $N_E$ is the Euclidean normal (identifying $H$ with
  $\R^3$) and $\mathscr{H}^2_E$ is the $2$-dimensional Hausdorff
  measure in $\R^3$.
\item In particular, if $S$ is given as a graph over a region $\Omega$
  of $u \in C^1(\Omega)$ and $\partial
  \Omega$ is also of class $C^1$ then
\[\mathscr{P}(u) = \int_\Omega |N_0| \;dx dy \]
\end{enumerate}
\end{Pro}
 
Franchi, Serapioni and Serra Cassano also prove an extension
of the implicit function theorem in the CC setting and use it
to show a beautiful structure theorem for $X$-Caccioppoli sets in $H$.
Again, for 
the purposes of this paper, we do not need to full power of their
result and will state only what is necessary.  To state the theorem,
we first recall some of the definitions in \cite{FSSC2}.  The reader
should be aware that we choose slightly different notation than that
in \cite{FSSC2} to maintain consistency within this paper.

\begin{Def} We call $S \subset H$ an {\bf $H$-regular hypersurface}
if, for every $p \in S$, there exists an open ball, $B_{cc}(p,r)$ and
a function $f:B_{cc}(p,r) \ra \R$ such that
\begin{enumerate}
\item Both $f$ and $\nabla_0 f$ are continuous functions.
\item $S$ is a level set of $f$, \[S \cap B_{c}(p,r) = \{q \in B_{cc}(p,r)
  | f(q)=0\}\]
\item $\nabla_0 f(p) \neq 0$
\end{enumerate}
\end{Def}
Roughly, $H$-regular surfaces are regular in the usual sense in the
distributional directions, but are allowed to be nondifferentiable in
the $Z$ direction.  

To define the reduced boundary, we need to define a generalized unit
normal to a surface base on the perimeter measure.  Given $E$, an
$X$-Caccioppoli set, consider the functional on $C^0_0$ given by 
\[ \phi \ra -\int_E div_0 \phi\; dV \]
where $div_0$ is the horizontal divergence operator, $(X,Y)$.  The
Riesz representation theorem implies that there exists a section of $\mathcal{V}$, $\nu_E$ so that

\[ -\int_E div_0 \phi \; dV = \int_H <\nu_E,\phi> d \mathscr{P}_E \]
$\nu_E$ is the {\bf generalized unit horizontal normal} to $\partial
E$.  In the case where $\partial E$ is smooth, $\nu_E$ coincides with $n_0$.

\begin{Def}  The {\bf reduced boundary} of an $X$-Caccioppoli set $E
  \subset H$, $\partial^*_{cc}E$ is the set of points $p \in E$ such
  that
\begin{enumerate} 
\item \[\mathscr{P}(E \cap B_{cc}(p,r)) >0\]
for all $r >0$.
\item \[\nu_0(p) = \lim_{r \ra 0} \frac{\int_{B_{cc}(p,r)} \nu_0 \; d
    \mathscr{P}_E}{\int_{B_{cc}(p,r)} \; d \mathscr{P}_E}\]
where $\mathscr{P}_E(\cdot) = \mathscr{P}(E \cap \cdot)$.
\item $|\nu_0(p)| =1$
\end{enumerate}
\end{Def} 
Note that the reduced boundary does not include any characteristic
points of the boundary.  

Now we state the pieces of the main theorems that we need for our
purposes. 
\begin{Thm}\label{FSSC}(\cite{FSSC2}, see theorems 6.4 and 7.1) If $E
  \subset H$ is an $X$-Caccioppoli set then,
\begin{enumerate}
\item \[\mathscr{P}_E = k \mathscr{H}^3_{cc} \hook \partial^*_{cc}E
  \] where $k$ is a constant.  
\item  $\partial^*_{cc}E$ is $H$-rectifiable, i.e.
\[ \partial^*_{cc}E = N \cup  \bigcup_{k=1}^\infty C_k\]
where $\mathscr{H}^3_{cc}(N)=0$ and $C_k$ is a compact subset on an
$H$-regular hypersurface $S_k$.
\end{enumerate}
\end{Thm}

\section{A Characterization of Minimal Surfaces}\label{char}
For the rest of the paper, we will restrict ourselves to considering
the three dimensional Heisenberg group equipped with the
Carnot-Carath\'eodory metric described in section \ref{back},
$(H,\cc)$.  Moreover, we will consider only nonparametric surfaces in $H$,
i.e. graphs of functions $u:\R^2 \ra \R$ over the $xy$-plane in $H$
(using the identification of $H$ with $\R^3$).  

\subsection{The Variational Setup}
To characterize minimal surfaces in $H$, we must first specify which
surface measure we wish to use.  As discussed in section \ref{back},
the perimeter measure, $\mathscr{P}$, is a natural measure to use in this setting and,
for sufficiently regular surfaces, it coincides with another natural
surface measure, the three dimensional spherical Hausdorff measure
restricted to the surface.  Proposition \ref{measform} yields that if
$S$ is the graph of $u: \Omega \subset \R^2 \ra \R$ where $u \in C^1(\Omega)$,
\[\mathscr{P}(u) = \int_\Omega |N_0| dx dy \]
where $N_0$ is the horizontal normal of $S$.  Taking this as our
starting point for a variational analysis, we define our energy
Lagrangian function  as 
\[ L(\xi,z,p) = ((\xi_1-y)^2+(\xi_2+x)^2)^\frac{1}{2}\]
Thus, for $u \in C^1(\Omega)$, 
\begin{equation*}
\begin{split}
L(\nabla u,u,\vec{x})&= ((u_x-y)^2+(u_y+x)^2)^\frac{1}{2}\\
&= ((XF)^2+(YF)^2)^\frac{1}{2} = |N_0|
\end{split}
\end{equation*}
where $F(x,y,z)=u(x,y)-z$.  Our energy function based on this
Lagrangian is
\[E(u) = \int_\Omega L(\nabla u, u, \vec{x}) \;dx dy \]
for $u \in W^{1,1}$.  Note that if $u \in C^1(\Omega)$ then $E(u)=\mathscr{P}(u)$.  For this energy function, the associated
Euler-Lagrange equation is
\[ - L_{\xi_1}(\nabla u,u,\vec{x})_x - L_{\xi_2}(\nabla u, u,
\vec{x})_y = 0 \]

Rewriting this in terms of the horizontal gradient operator, $\nabla_0$, this
yields:
\begin{equation*}\label{MSE}
-\nabla_0 \cdot \frac{\nabla_0 F}{|\nabla_0 F|} =0 \tag{MSE}
\end{equation*}

Note that this is simply the horizontal 1-Laplacian on the Heisenberg group.
In keeping with classical notation, we define a partial differential
operator:
\[ H_{cc}(u) = \nabla_0 \cdot \frac{\nabla_0 F}{|\nabla_0 F|} \]

\begin{Lem}$L$ is convex in $\xi$.
\end{Lem}
\pf The Hessian of $L$ is positive semidefinite with eigenvalues $\{0,
((\xi_1-y)^2 + (\xi_2+x)^2)^{-\frac{1}{2}}\}$. $\qed$

By the
standard variational theory we have that $E(\cdot)$ is weakly lower
semicontinuous on $W^{1,p}$.  Moreover, the convexity implies
that any weak solution to $(\ref{MSE})$ is a minimizer of the energy
function $E(\cdot)$.  A priori, these ``minimal
surfaces'' may not be the $X$-minimal surfaces of Garofalo and Nhieu
(theorem \ref{GN}) described in section \ref{back}.  However, we shall see that
the solutions to the Plateau problem in the next section are indeed
$X$-minimal surfaces.
We note that it is possible, albeit cumbersome, to give a 
classical derivation of the same characterization of minimal surfaces
assuming that the surfaces are at least $C^2$ and are critical with
respect to compactly supported $C^2$ variations.

We point out that because the equation $(\ref{MSE})$ is not strictly
elliptic, we do not automatically get uniqueness of solutions with
prescribed boundary data.  Indeed, as we shall see in the subsequent
sections, solutions are not necessarily unique.

\subsection{Minimal surfaces as the limit of approximating minimal
  surfaces}

In this section, we examine the connection between minimal surfaces in
the approximating spaces $(H,d_\lambda)$ and the minimal surfaces in
$(H,\cc)$.  First, we illustrate that some minimal surfaces in
$(H,\cc)$ arise as limits of sequences of minimal surfaces, each of
which is in an $(H, d_\lambda)$.  Second, we use this description to
give another (more geometric) proof of the existence of solutions to
the Plateau problem in $(H,\cc)$ by constructing nonparametric
solutions via limits of solutions to Plateau problems in the
approximating spaces.    

In $(H,d_\lambda)$, we use the standard Lagrangian and energy
functions used to investigate minimal surfaces:
\begin{equation*}
L_\lambda(\xi,z,p)=
\left((\xi_1-y)^2+(\xi_2+x)^2+\frac{1}{\lambda^2}\right)^\frac{1}{2}
\end{equation*}
\[ E_\lambda(u)=\int_\Omega L_\lambda(\nabla u, u,\vec{x})dx dy \]

Later, we will also compare this functional with the Riemannian area
functional on surfaces in $(H,g_\lambda)$.  To this end, we define,
for a smooth surface $S$ in $(H,g_\lambda)$,
\[A_\lambda(S) = \int_S dA_\lambda \]
where $dA_\lambda$
is the Riemannian area element.  Note that if $S$ is given by the
graph of a function $u:\Omega \ra\R$ then
$A_\lambda(u)=E_\lambda(u)$.  

\begin{Lem}\label{rel1}  Given a function $F(x,y,z)=u(x,y)-z$ on $H$ where
  $u:\Omega \subset \R^2 \ra \R$, 
\[E(u) \le E_\lambda(u) \le E(u)+\frac{1}{\lambda} \int_\Omega dxdy \]
\end{Lem}
\pf This follows from the definition of the energy functions since
$L_\lambda=((\xi_1-y)^2+(\xi_2+x)^2+\frac{1}{\lambda^2})^{\frac{1}{2}}$ and
$L=((\xi_1-y)^2+(\xi_2+x)^2)^\frac{1}{2}$. $\qed$

\begin{Pro}\label{lims}  Consider the nonparametric minimal surface
  equation in $(H, d_\lambda)$, $H_\lambda =0$.  Suppose for each
  $\lambda >1$, there exists $u_\lambda: \Omega \subset \R^2
  \ra R$ so that $u_\lambda$ minimizes $E_\lambda(\cdot)$.  Further
  suppose that $\{u_\lambda\}$ converges weakly to some function $u$ in
  $W^{1,p}(\Omega)$.  Then $u$ is a minimizer of $E(\cdot)$ and
  $F=u-z$ is
  a weak solution to $(\ref{MSE})$.
\end{Pro}

\pf  Let $m=\inf_g E(g)$ and $m_\lambda = \inf_g E_\lambda(g)$.  Since
each $u_\lambda$ is a minimizer of $E_\lambda$, by the previous lemma,
we have 
\begin{equation*}
E(u_\lambda) \le m_\lambda \le E(u_\lambda)+ \frac{1}{\lambda}
\int_{\Omega} dx dy
\end{equation*}
Since $m \le E(u_\lambda)$ and, if $u$ is a minimizer of $E(\cdot)$
with the same domain,
\begin{equation*}
m \le E_\lambda(u) \le m +\frac{1}{\lambda} \int_\Omega \;dx dy 
\end{equation*}
Combining these, we have
\[ m \le m_\lambda \le m+ \frac{C}{\lambda}\]
where $C$ is a constant depending on $\Omega$.  Thus, as $\lambda \ra \infty$, we have \[ \lim_{\lambda \ra \infty}
E_\lambda(u_\lambda)=m \]
In other words, using this and lower semicontinuity, $\{u_\lambda\}$ converges to a minimizer of
$E(\cdot)$.  Hence, $F=u-z$ satisfies $(\ref{MSE})$ weakly.  $\qed$

Using standard elliptic estimates, we can produce a weak
solution to the Plateau problem in $(H,\cc)$ for curves satisfying the
bounded slope condition.  We first recall the bounded slope condition
in this setting:

\begin{Def}  Let $\Gamma$ be a closed curve in $H$ given as the graph
  of a function $\phi:\partial \Omega \ra \R$ where $\Omega$ is a
  domain in $\R^2$.  Then $\Gamma$ has the {\bf bounded slope
  condition} if, for every point $p \in B=\{ (x,y,z) \in H | (x,y) \in
  \partial \Omega, z=\phi(x,y)\}$, there exist two planes $P_{-}(p)$ and
  $P_+(p)$ passing through $p$ such that
\begin{enumerate}
\item \[P_-(p) \le \phi(x,y) \le P_+(p)\]
\item Let $s(P)$ denote the slope of a plane.  Then, the collection
  \[\{s(P_-(p))|p\in B\}\cup\{s(P_+(p))|p\in B\}\]
is uniformly bounded by some constant $K$.
\end{enumerate}
\end{Def}
 
\begin{Thm}\label{minlims}  Let $\Gamma$ be a closed curve in $H$ satisfying
  the bounded slope condition and which is the graph
  of a function $\phi \in C^{2,\alpha}(\R)$ over a curve $\gamma \in C^{2,\alpha}(\R^2)$ which bounds a
  region $\Omega$.  Then, there exists $u \in W^{1,p}(\Omega) \cap C^0(\Omega)$ so that
  $u|_\gamma = \phi$, $u$ is a weak solution to $(\ref{MSE})$ on
  $\Omega$ and $u$ minimizes $E(\cdot)$ on $\Omega$.  Moreover, there
  exists a sequence of functions with the same regularity as $\Gamma$,
  $\{u_{\lambda_n}\}$, such that $H_{\lambda_n}(u_{\lambda_n})=0$ and
  $u_{\lambda_n} \ra u$ in $W^{1,p}(\Omega)$.
\end{Thm}

\pf For a fixed $\lambda \ge 1$, $H_\lambda =0$ is a uniformly elliptic
equation and so, by the standard theory (see, for example
\cite{GilTrud}, theorem 11.5), since $\Gamma$ is a $C^2$ curve satisfying the
bounded slope condition, there exists a
solution $u_\lambda$ to the Plateau problem in $(H,d_\lambda)$.  In
other words, $H_\lambda(u_\lambda)=0$ on $\Omega$, $u_\lambda|_\gamma =
\Gamma$ and $u_\lambda$ minimizes $E_\lambda(\cdot)$.  Moreover, standard applications of the quasilinear elliptic maximum principle yield:
\begin{equation*}
\begin{split}
\sup_\Omega |u_\lambda| & \le \sup_\gamma |u_\lambda| = \sup_\gamma
|\phi|\\
\sup_\Omega |\nabla u_\lambda| & \le \sup_\gamma |u_\lambda| =
\sup_\gamma |\phi|
\end{split}
\end{equation*}
Thus, each $u_\lambda$ is bounded in $W^{1,p}$ by a constant $C$
depending only on $\phi$.  In particular, $C$ does {\it not} depend on
$\lambda$.  Thus, $\{u_\lambda\}$ is uniformly bounded in
$W^{1,p}(\Omega)$ and hence, we can extract a subsequence
converging weakly to a limit function, $u_\infty \in W^{1,p}_0$.
In fact, since the $u_\lambda$ are all $C^2$, the same boundedness
yields that the family is uniformly bounded and equicontinuous and
hence we may assure that the convergence is uniform.  Thus, $u_\infty$ is
at least continuous, $u_\infty|_\gamma=\phi$ and
$u_\infty$ is a graph over $\Omega$.  By the previous proposition,
$F=u_\infty(x,y)-z$ is a weak solution to the Plateau problem in
$(H,\cc)$ with the specified boundary data. $\qed$\\

The previous theorem yields a nice characterization of some solutions
to the Plateau problem - they are limits of solutions to the Plateau
problem in Riemannian approximates to $(H,\cc)$.  However, we now show
a modification that shows that one can generate solutions with much
more flexibility. 

\begin{Thm}\label{meanlims}   Let $\Gamma$ be a closed curve with a
bounded slope condition in $H$ which is the graph
  of a function $\phi\in C^{2,\alpha}(\R)$ over a curve $\gamma \in C^{2,\alpha}(\R^2)$ which bounds a
  region $\Omega$.  Suppose $u_n$ are graphs in $(H,d_{\lambda_n})$
  spanning $\Gamma$ which have mean curvature given by the function
  $\kappa_n$.  If $\kappa_n \ra 0$ as $n \ra \infty$
  and $u_n$ converge to a graph, $u$, then $u$ satisfies (\ref{MSE}) and
  is a minimizer of $E(\cdot)$.  
\end{Thm}

\pf First, since $u_n$ is a
solution to a prescribed mean curvature equation with boundary values
given by $\Gamma$, its energy $E_{\lambda_n}(u_n)$ is close to
$E(u_n)$ because the Lagrangian formulation of the prescribed mean
curvature problem uses $E_\lambda(\cdot)$ as the energy coupled with
an additional integral constraint.  For a fixed $\lambda_n$, if a
convergent sequence of surfaces have mean curvatures converging to zero
(uniformly), then they converge to a minimal surface.  Moreover, since
solutions to the minimal surface Dirichlet problem in $(H,g_{\lambda_n})$ exist (again due to
the bounded slope condition of $\Gamma$) and are unique (see theorems
10.1 and 10.2 in \cite{GilTrud}), they
must converge to the unique energy minimizer.  Thus, 
\[ |E_{\lambda_n}(u_n) - E(u_n)| = o_{n^{-1}}(1) \]

Second,  we
follow the argument in proposition \ref{lims}.  Letting $u_{\lambda_n}$ be
the minimal surface spanning $\Gamma$ in $(H,d_{\lambda_n})$, $u$ a
minimizer of $E(\cdot)$ spanning $\Gamma$,
$m_n=E_{\lambda_n}(u_{\lambda_n})$ and $m=\inf_g E(g)$, we see that by lemma \ref{rel1}
\begin{equation*}
\begin{split}
E_{\lambda_n}(u_n) &= m_n + o_{n^{-1}}(1)\\
& \ge E(u_{\lambda_n})+ o_{n^{-1}}(1)\\
& \ge m +  o_{n^{-1}}(1)
\end{split}
\end{equation*}
Moreover, again using lemma \ref{rel1} and the argument in proposition
\ref{lims}, we have
\begin{equation*}
\begin{split}
E_{\lambda_n(u_n)} &= m_n + o_{n^{-1}}(1)\\
&\le E_{\lambda_n}(u) + o_{n^{-1}}(1)\\
&\le m+\frac{C}{\lambda_n} + 
 o_{n^{-1}}(1)
\end{split}
\end{equation*}
Taking these together, we have
\[ m+ o_{n^{-1}}(1) \le E_{\lambda_n}(u_n) \le m+ \frac{C}{\lambda_n}+
o_{n^{-1}}(1)\]

Hence, as $n \ra \infty$, $\lambda_n \ra \infty$, and $E_\lambda(u_n)$ converges to
$m=E(u)$.  Thus, by lower semicontinuity, we have that $E(u_n)=m$ and
hence is a minimizer of $E(\cdot)$ as well.\\
$\qed$

\begin{Rem}  We point out one feature of the partial differential
  operator $H_{cc}$, namely that if $H_{cc}(u)=0$, then  
\[H_\lambda(u) = \frac{1}{\lambda^2}
\frac{\Delta u}{\left((u_x-y)^2+(u_y+x)^2+\frac{1}{\lambda^2}\right
  )^\frac{3}{2}}\] where $\Delta$ is the usual Laplacian on  $\R^3$.
  Therefore, any solution to $H_{cc}=0$ with boundary data satisfying
  the bounded slope condition satisfies the hypotheses of theorem
  \ref{meanlims}.  As we generate examples of solutions to $H_{cc}=0$ in
  section \ref{egs}, this observation
  provides a link to the energy minimizers in this section.
\end{Rem}

Next, we provide a link to the $X$-minimal surfaces given
in \cite{GarNh}.  
\begin{Thm}\label{min2} The minimizers of $E(\cdot)$ constructed in
  theorems \ref{minlims} and \ref{meanlims} are also
  X-minimal surfaces.  In other words, they are also minimizers for
  the perimeter measure.
\end{Thm}

\pf  By the lower semicontinuity of the perimeter measure and the
characterization of the perimeter measure for smooth functions in
lemma \ref{semicont}, we see that the perimeter of an open set, $O$, with boundary given
by $u$ (over $\Omega$) is larger or equal to the perimeter of one of Garofalo
and Nhieu's $X$-minimal surfaces with the same boundary conditions.
Conversely, given $W$ an $X$-minimal open set with boundary conditions matching that of $u$, we
can use the work of Franchi, Serapioni and Serra Cassano (\cite{FSSC})
as follows.  By theorem \ref{FSSC}, we can realize the reduced
boundary of $W$ by level sets of $H$-regular functions up to a set of
$\mathscr{H}^3_{cc}$ measure zero.  Without loss of generality, we can
focus only on one of the level sets given by a function $\psi:\Omega
\subset \R^2 \ra
\R$ - we may have to piece together many of these patches, but the
result will be the same up to a set of measure zero.  Moreover, using
proposition \ref{smoothapprox}, we can produce a smooth
approximate $\psi_\epsilon$ of $\psi$ which converges in $L^1$ to
$\psi$ as $\epsilon \ra 0$ with the property that \[\lim_{\epsilon \ra
  0} \mathscr{P}(\psi_\epsilon) = \mathscr{P}(\psi)\]

Unfortunately, $\psi_\epsilon$ may no
longer have the same boundary data as $\psi$.  To adjust this, let
$T_\epsilon$ be an $\epsilon$ neighborhood of $\partial \Omega$ and
construct a smooth function $h_\epsilon$ with support in $T_\epsilon$ and so
that $\psi_\epsilon-(h_\epsilon-z)=0$ restricted to $\partial \Omega$
coincides with the graph of $\phi$.  Note that we may arrange that
$h_\epsilon-z$ has bounded gradient with a bound not depending on
$\epsilon$ and that $h_\epsilon 
\ra 0$ as $\epsilon \ra 0$.  Notice that the perimeter of
$\psi_\epsilon$ and $\psi_\epsilon-(h_\epsilon-z)$ can be computed using
$\mathscr{P}$ by proposition \ref{measform} as both functions are
smooth.  Further, letting $W_\epsilon$ be the open set
$\psi_\epsilon<0$ and $W_{h,\epsilon}$ be the open set
$\psi_\epsilon-(h_\epsilon-z)<0$, 

\begin{equation*}
\begin{split}
|\mathscr{P}(W_{h,\epsilon})-\mathscr{P}(W_\epsilon)|
&\le \int_{T_\epsilon \cap \Omega} ||\nabla_0(\psi_\epsilon -
(h_\epsilon-z))| -|\nabla_0 (\psi_\epsilon)||\; dx dy\\
&\le \int_{T_\epsilon \cap \Omega} |\nabla_0 (h_\epsilon-z)| \;dx dy\\
&\le C \int_{T_\epsilon \cap \Omega}\; dx dy
\end{split}
\end{equation*}

Thus, 
\begin{equation}\label{a1}
|\mathscr{P}(W_{h,\epsilon})-\mathscr{P}(W_\epsilon)|
\ra 0 \text{\;\;as $\epsilon \ra 0$}
\end{equation}
Lastly, the relation between $A_\lambda$ and $E_\lambda$ coupled with the
representation of $\mathscr{P}$ for smooth functions given in
proposition \ref{measform} yields:
\begin{equation}\label{a3}
|A_\lambda(S)-\mathscr{P}(O)| = o_{\lambda^{-1}}(1)
\end{equation}
for any smooth surface $S$ bounding an open set $O$.

Now, since $\partial W_{h,\epsilon}$ is smooth,
using proposition \ref{measform}, we see that 
\begin{equation*}
\begin{split}
\mathscr{P}(W_\epsilon) &= \mathscr{P}(W_{h,\epsilon})
+o_{\epsilon}(1) \text{\;\; by (\ref{a1})}\\
&=
A_\lambda(\psi_\epsilon-h_\epsilon)+o_{\lambda^{-1}}(1)+o_{\epsilon}(1)
\text{\;\; by (\ref{a3})}\\ 
&\ge
A_\lambda(u_\lambda)+o_{\lambda^{-1}}(1)+o_{\epsilon}(1)\\
&=E_\lambda(u_\lambda)+o_{\lambda^{-1}}(1)+o_{\epsilon}(1)
\end{split}
\end{equation*}
where $u_\lambda$ is a
Riemannian minimal graph with the same boundary data minimizing
$E_\lambda(\cdot)$.  By the argument in the proof of proposition
\ref{lims} and proposition \ref{measform}, we have that
\[E_\lambda(u_\lambda) \ge E(u_\lambda) = \mathscr{P}(u_\lambda)\]
By the lower semicontinuity of $\mathscr{P}$ (lemma \ref{semicont}),
we have \[\liminf_{\lambda \ra \infty}
\mathscr{P}(u_\lambda) \ge \mathscr{P}(u)\]
Putting this together with the previous calculation we have, finally,
\[\mathscr{P}(W) = \lim_{\epsilon \ra 0} \mathscr{P}(W_\epsilon) \ge
\mathscr{P}(u)\]  Thus, $u$ and $W$ describe sets of minimal
perimeter.  $\qed$ \\

{\em Remarks:}
\begin{enumerate}
\item These minimal surfaces in $(H,\cc)$ are geometrically realizable
  as limits and can often be explicitly constructed.  In section
  \ref{egs}, many of the examples can be realized in this way.
\item Theorem $\ref{meanlims}$ shows that examples of minimal surfaces
  in $(H,\cc)$ can arise in a wide variety of limits.  This, along with the
subellipticity of $H_{cc}=0$ is strong evidence that the solution to
the Plateau problem is not neccesarily unique.  
\item Of course, one would like the solutions to have higher
  regularity than given in the theorem.  However, once again the
  examples in section \ref{egs} (particularly those in section
  \ref{Gauss}) show that, without imposing boundary conditions,
  regularity cannot be expected.  The question of higher regularity of
  solution to the Plateau problem with $C^k$ boundary data is open.
\end{enumerate}

\section{Examples}\label{egs}
\subsection{Invariant Solutions}\label{inv}
In \cite{Tomter}, P. Tomter began a program of finding minimal
surfaces in $(H,g_1)$ which are invariant under isometric group
actions based on the general technique of Hsiang and Lawson
(\cite{Hsiang:Lawson}), later expanded by Hsiang and Hsiang (\cite{Hsiang:Hsiang}).  This program is completed in \cite{FMP} by Figuero, et
al. where a number of additional cases are studied.  The method rests
on reducing the minimal surface partial differential equation to an
ordinary differential equation (or system of ODEs) by considering how
the PDE descends to the quotient of $H$ by a closed subgroup of
isometries.  Via the Riemannian submersion given by the process of
quotienting by a group of isometries, the solutions of the ODE(s) can
be lifted to minimal surfaces in $(H,g_1)$.  For complete details on
the method in the Heisenberg group, see \cite{Tomter} and \cite{FMP}.  

More specifically, the authors mentioned above consider subgroups of
rotations, group translations and combinations of these, which we will
refer to as ``corkscrew'' motions.  In light of proposition
\ref{lims} and the fact that the subgroups of rotations, translations
and corkscrew motions are also subgroups of the isometry group of
$(H,\cc)$, one expects that there exist group invariant minimal
surfaces in $(H,\cc)$.  In this section, we identify, by brute
force, several families of such minimal surfaces.  We should remark
that our technique is not a generalization of the quotienting
technique alluded to above - the quotient mapping is a type of
submersion, but it does not behave nearly as well as the Riemannian 
version.  Thus, for the purposes of this paper, we rely on the
following two methods:

\begin{enumerate}
\item Search for examples as limits of Riemannian examples using
  proposition \ref{lims} or theorem \ref{meanlims}.
\item Search for invariant examples directly in $(H,\cc)$.
\end{enumerate}

While the first method does reveal specific initial examples, the
second is much more efficient in finding the most examples.   We will
first give some simple isolated examples of the first technique and then
explore the second.  
\subsubsection{Invariant solutions as limits}

\vspace{.1in}

\noindent
{\em Planes and Helicoids.}
Consider the nonparametric minimal surface equation in
$(H,g_\lambda)$, $H_\lambda(u)=0$ or:
\begin{gather}
u_{xx}\left ( \frac{1}{\lambda^2}+(u_y+x)^2 \right )
-2u_{xy}(u_y+x)(u_x-y) +u_{yy}\left ( \frac{1}{\lambda^2} + (u_x-y)^2
\right )=0  \label{MSL}
\end{gather}

By inspection, we see that any plane $z=ax+by+c$ is a solution to
(\ref{MSL}) for any $\lambda >0$.  Thus, by proposition \ref{lims},
all planes are solutions of (\ref{MSE}).  

Further, considering $\R^3$ using cylindrical coordinates and defining
the helicoids by $u(r,\theta)=a\theta$ for $a \in \R$, then it is a quick computation
to show that all of the helicoids satisfy (\ref{MSL}).

%{\bf A nontrivial limit example?}

\subsubsection{Direct methods}
Under the assumption that invariant minimal surfaces exist, in this
section we produce smooth nonparametric examples.

\vspace{.1in}

\noindent
{\em Rotationally and corkscrew invariant examples:}  Again viewing $\R^3$ in cylindrical coordinates, we look for solutions
of the form $u(r,\theta)=f(r)$.  Under the transformation to
cylindrical coordinates, (\ref{MSE}) becomes:

\begin{equation}
\frac{u_{rr}(r^5+ru_\theta^2+2r^3u_\theta)-u_{\theta
    r}(2ru_\theta u_r+2r^3u_r)+u_{\theta
    \theta}u_r^2+2u_ru_\theta^2+2r^2u_\theta u_r+r^2u_r^3}{r^3(u_r^2+r^{-2}u_\theta^2+2u_\theta+r^2)^{\frac{3}{2}}}
    =0  \label{Mcyl}
\end{equation}

We are looking for all graphs invariant under a combination of rotations
and vertical translations (the corkscrew motions), namely
\[u(r,\theta) = v(r)+a \theta \text{\;\; $a \in \R$}\]
Under this assumption, (\ref{Mcyl}) becomes:
\begin{equation}
\frac{v''(r)(ra^2+2r^3a+r^5)+v'(r)(2a^2+2r^2a)+r^2v'(r)^3}{r^3(v'(r)^2+\frac{a^2}{r^2}
  +2a+r^2)^\frac{3}{2}}=0
\end{equation}
or, under the assumption that the denominator does not vanish,
\begin{equation}\label{cork1}
v''(r)= -\frac{v'(r)(2a^2+2r^2a)+r^2v'(r)^3}{r(a^2+2r^2a+r^4)}
\end{equation}

This ordinary differential equation can be solved in closed form
yielding
\begin{equation}\label{corksol}
v(r)=\pm \left (\frac{\sqrt{br^2-1}}{b}-a\tan^{-1}\left (
  \frac{1}{\sqrt{br^2-1}}\right )\right )+c \text{\; \; for $b,c, \in
  \R$}
\end{equation}

Writing the graphs implicitly yields:

\[z^2+2a\theta z=\frac{r^2}{b} - \frac{1}{b^2} - \frac{2a}{b}\tan^{-1}\left (
  \frac{1}{\sqrt{br^2-1}}\right )\sqrt{br^2-1}+a^2\left (\tan^{-1}\left (
  \frac{1}{\sqrt{br^2-1}}\right )\right )^2 - a^2\theta^2\]

\begin{figure}
\centering 
\mbox{\subfigure[$a=0,b=1$]{\epsfig{figure=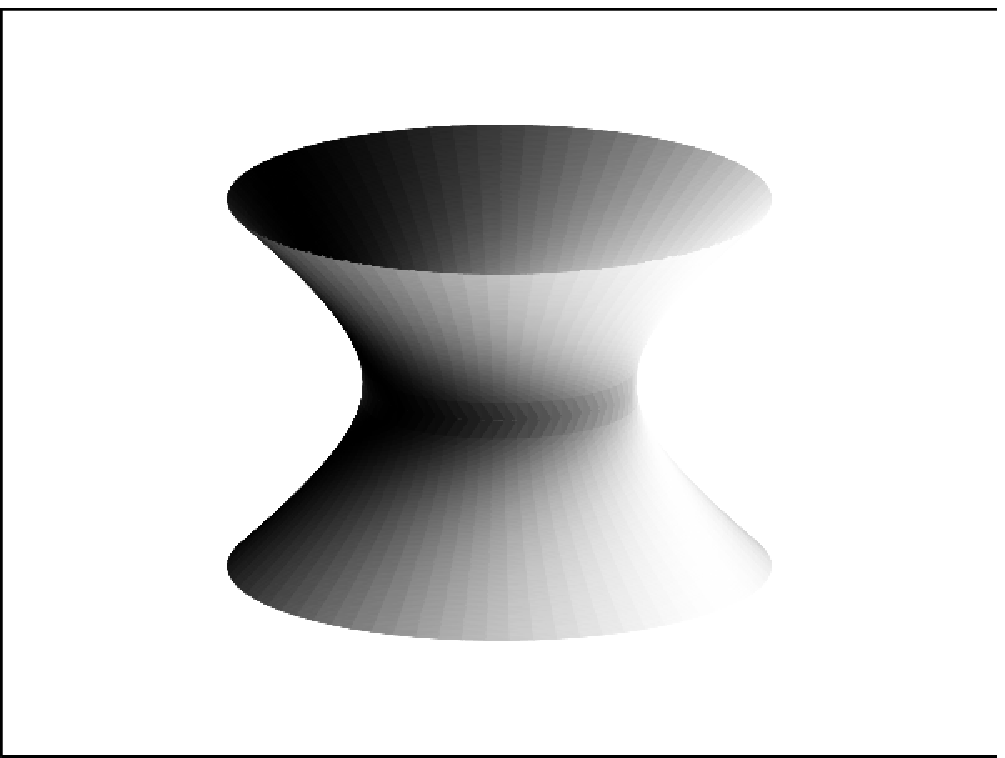,width=2.7in}}\quad
\subfigure[$a=1,b=1$]{\epsfig{figure=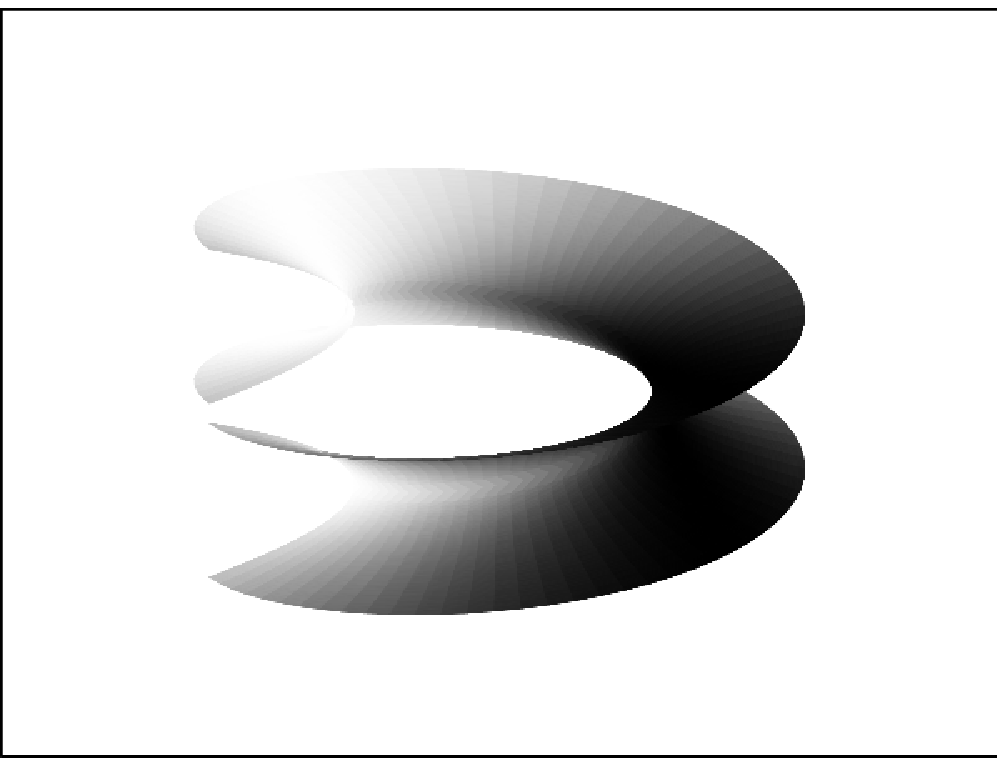,width=2.7in}}}
%\mbox{
%\subfigure[$a=-1,b=1$]{\epsfig{figure=cork2.eps,width=2.7in}}}
\caption{Corkscrew invariant minimal surfaces}\label{corkpics}
\end{figure}

Figure \ref{corkpics} shows two examples of this type.  Notice
that when $a \neq 0$, these surfaces are not embedded, but immersed.
The shading in the figures shows the change of the horizontal normal -
the point is darker the closer  $<n_0,X>$ is to $1$.

\vspace{.1in}

\noindent
{\em Translationally invariant examples:}  

We now consider graphs which are invariant under left translations by
one parameter subgroups of the Heisenberg group.  There are two cases,
when we translate in the vertical direction and when we translate in
a horizontal direction.  Since a surface invariant under a vertical
translation would not be a graph, we will focus on the second case.
Due to the rotational invariance of the metric, it suffices to
consider only graphs invariant under the subgroup of isometries
generated by translation in the $X$ direction.  Thus, we look for
graphs $z=u(x,y)$ so that $(x+t,y,u(x+t,y))=(t,0,0)\cdot (x,y,u(x,y))$.  This
amounts to solving the equation:
\[u(x+t,y)= u(x,y)+ty\]

One quickly sees that any smooth $u$ is of the form $u(x,y)=xy+f(y)$
where $f$ is any smooth function.  Under this assumption, the minimal
surface equation (\ref{MSE}) is identically zero.  Thus, any graph of
the form \[z=xy+f(y)\] satisfies the minimal surface equation.  

\noindent
{\em Remarks}: 
\begin{enumerate}
\item Note that so long as $f$ has no singularities, all of
  these examples are {\em complete} minimal graphs.
\item The horizontal unit normal of these surfaces is a piecewise
  constant vector.  Indeed, $X(xy+f(y)-z)=0$ and
  $Y(xy+f(y)-z)=2x+f'(y)$.  So, the horizontal unit normal,
  \[n_0=\left (0,\frac{2x+f'(y)}{|2x+f'(y)|}\right)\]  Notice that the vector field
  is discontinuous along the characteristic locus over which it changes
  sign. 

\end{enumerate}

\subsection{Prescribed Gauss maps}\label{Gauss}
In classical minimal surface theory, the Gauss map plays a crucial
role in recognizing and constructing examples of minimal surfaces.  In
this section, we analyze the horizontal Gauss map to help produce
families of examples of minimal surfaces derived from the initial
examples generated in the previous section.  

Recall the definition of the horizontal Gauss map for a nonparametric
surface, $F(x,y,z)=u(x,y)-z=0$:
\begin{equation}\label{N0}
n_0 = \frac{\nabla_0 F}{|\nabla_0 F|} = \left (
  \frac{p}{\sqrt{p^2+q^2}}, \frac{q}{\sqrt{p^2+q^2}} \right )
\end{equation}

Where $p = u_x-y$ and $q=u_y+x$.  Using this notation, $(\ref{MSE})$
becomes
\begin{equation}\label{divform}
div \; n_0 =0
\end{equation}
where {\em div} is the standard divergence operator in $\R^2$.  Observing
this, one strategy for finding examples is to search for divergence
free unit vector fields on $\R^2$ and associate to them solutions of
$(\ref{divform})$.  If $V=(v_1,v_2)$ is such a vector field, one can use
the following procedure to generate examples:
\begin{enumerate}
\item From $(\ref{N0})$ and $(\ref{divform})$, we see that we must have
\[v_1= \frac{p}{\sqrt{p^2+q^2}} \text{ \; \; and \; \;}
v_2=\frac{q}{\sqrt{p^2+q^2}}\]

\item Attempting to solve the equations above for $u(x,y)$, one
  possible solution is for $p$ and $q$ to satisfy $v_2p=v_1q$.
\item This gives rise to a first order partial differential equation:
\begin{equation}\label{linpde}
v_2 (u_x-y)-v_1(u_y+x)=0
\end{equation}
\end{enumerate}

Thus, any solution to $(\ref{linpde})$ yields a solution to $(\ref{divform})$
at points where $\nabla_0F \neq 0$.  At these, the characteristic
points of the graph, the Gauss map is discontinuous.  Indeed, a quick
calculation shows that if $v_2p=v_1q$,
\[n_0 = \left ( \frac{p}{|p|} v_1, \frac{q}{|q|} v_2 \right )\]

Thus, depending on the values of $p$ and $q$ when crossing over the
characteristic locus, the value of $n_0$ may flip, for example, from
$V$ to $-V$.  

In summary, we have,

\begin{Thm}If $V=(v_1,v_2)$ is a divergence free unit vector field on
  $\R^2$ and $u(x,y)$ is a solution to 
\begin{equation}
v_2(u_x-y)-v_1(u_y+x)=0 \tag{\ref{linpde}}
\end{equation}
Then $u(x,y)$ satisfies $(\ref{MSE})$ and its Gauss map, $n_0$,
matches $V$ up to sign except at characteristic points, where it is
discontinuous.
\end{Thm}

This theorem allows us to construct many surprising examples.

\subsubsection{Constant Gauss map}\label{const}

In the initial set of examples generated above, we noticed that of the
translationally invariant examples, the surfaces $z=xy+f(y)$ had
piecewise constant horizontal Gauss map.  We will now apply the
theorem above to find more examples with piecewise constant Gauss map.  As we
will discover, all of these examples are complete, in other words,
they are solutions to the Bernstein problem - complete minimal graphs.
To construct them, we simply prescribe 
\[ V=(\alpha, \pm \sqrt{1 - \alpha^2}), \; \alpha \in [-1,1]\]
Taking first the positive sign on the square root, $(\ref{linpde})$ becomes
\[   \sqrt{1 - \alpha^2}(u_x -y) - \alpha (u_y+x)=0\]
Using the method of characteristics and computing the envelope of the
resulting solution, we find the following solutions to the above
linear partial differential equation:\\
\begin{center}
\boxed{
u(x,y)=
\frac{\alpha}{\sqrt{1-\alpha^2}}x^2+xy+g\left
  (y+\frac{\alpha}{\sqrt{1-\alpha^2}}x\right )
}\\
\end{center}
where $g$ is {\em any} function in $W^{1,p}(\R)$.  In this special
case, we have
\begin{equation*}
p =
2\frac{\alpha}{\sqrt{1-\alpha^2}}x+\frac{\alpha}{\sqrt{1-\alpha^2}}Dg\left(y+\frac{\alpha}{\sqrt{1-\alpha^2}}x\right)
\end{equation*}
and
\begin{equation*}
q= 2x+Dg\left(y+\frac{\alpha}{\sqrt{1-\alpha^2}}x\right)
\end{equation*}
and so 
\[n_0 =
\frac{2x+Dg\left(y+\frac{\alpha}{\sqrt{1-\alpha^2}}x\right)}{|2x+Dg\left(y+\frac{\alpha}{\sqrt{1-\alpha^2}}x\right)|}
(\alpha, \sqrt{1-\alpha^2}) \]

\begin{figure}
\centering 
\mbox{\subfigure[$g(s)=0$]{\epsfig{figure=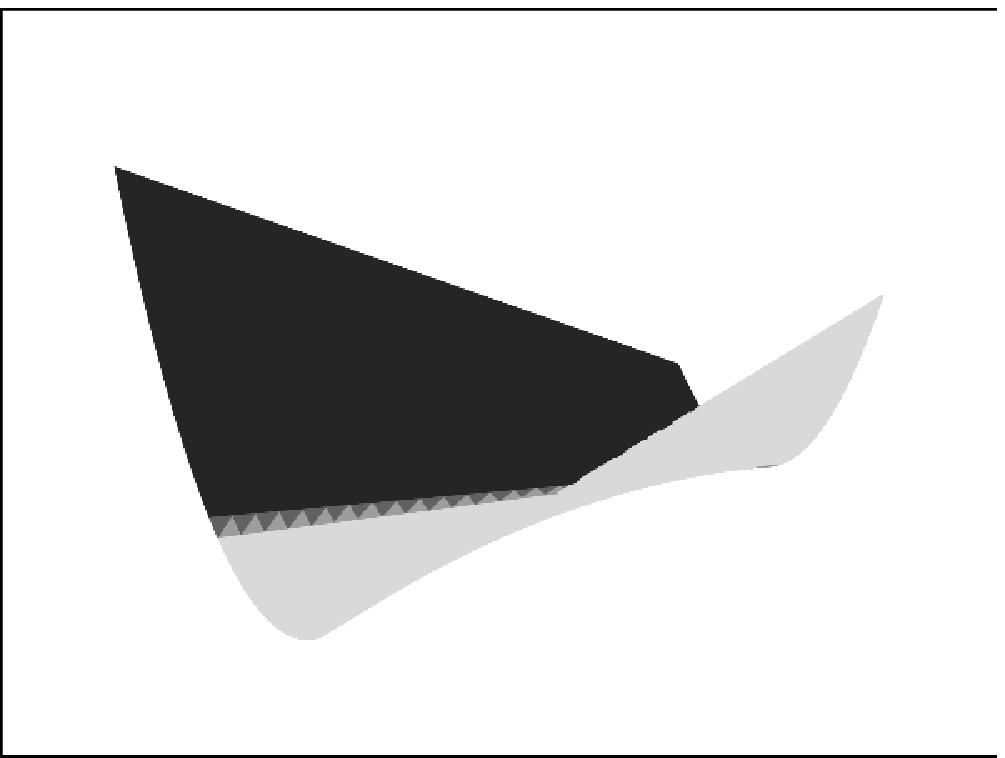,width=2.7in}}\quad
\subfigure[$g(s)=-s^2$]{\epsfig{figure=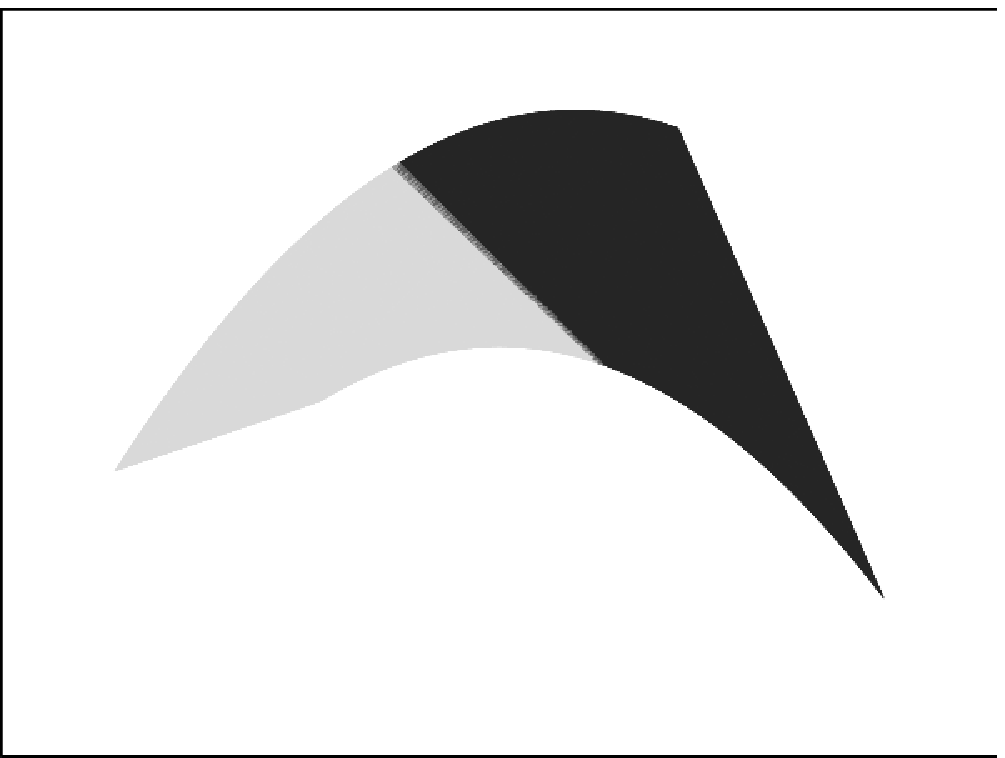,width=2.7in}}}
\mbox{
\subfigure[$g(s)=2\sin(3s)^2$]{\epsfig{figure=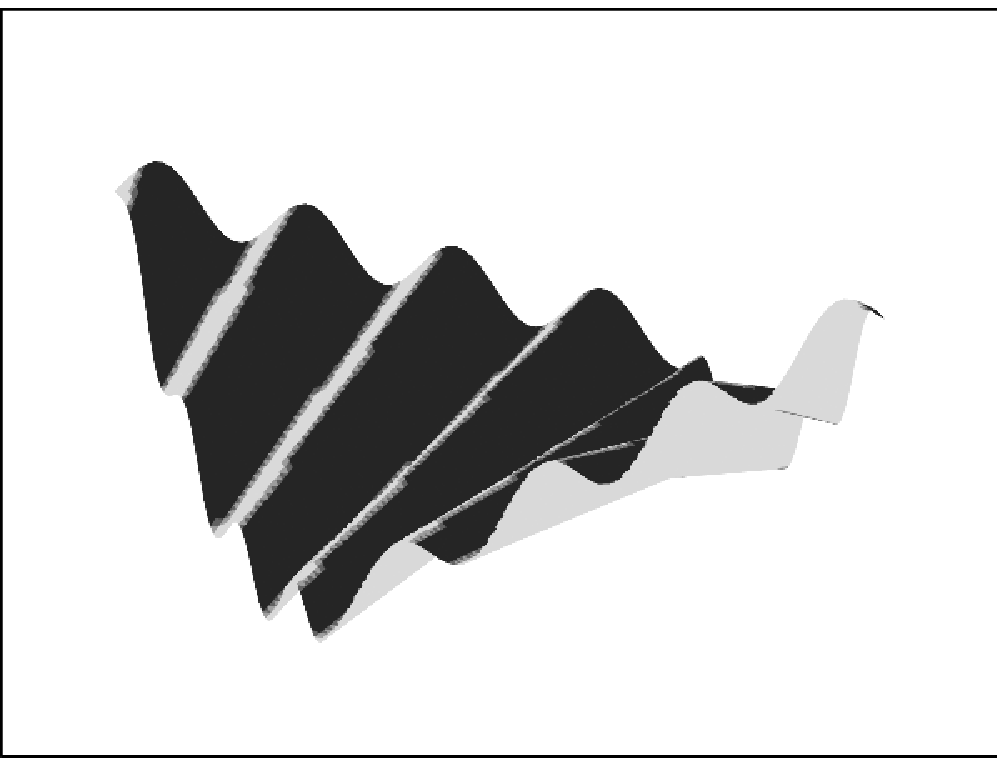,width=2.7in}}\quad
\subfigure[$g(s)=8e^{-\frac{(y+x)^2}{10}}$]{\epsfig{figure=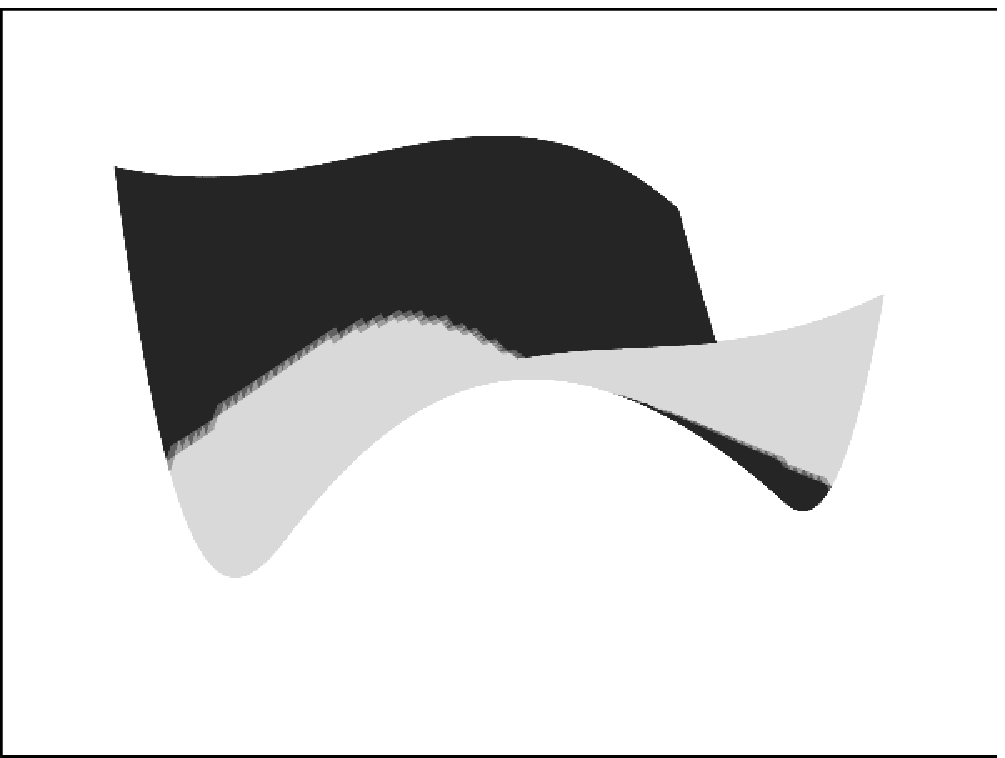,width=2.7in}}}
\caption{Minimal surfaces with constant Gauss map}\label{constN}
\end{figure}

Figure \ref{constN} show the graphs of several examples with different
choices of $g$ and with $\alpha=\frac{\sqrt{2}}{2}$.  The shading in
the figures represents the different patches of constant Gauss map
which differ from one another in sign.  The interface of the darker
and lighter shades is precisely the characteristic locus.  There is one degenerate case not covered above, we
$\alpha = \pm 1$.  In the case $\alpha=1$, $(\ref{linpde})$ becomes
$q=u_y+x=0$ which implies $u(x,y)=-xy+g(x)$.  So, $p=-2y+Dg(x)$ and
\[n_0 = \left (\frac{ -2y+Dg(x)}{|-2y+Dg(x)|},0\right)\]
The case $\alpha=-1$ yields the same solutions. 

If one takes $V=(\alpha, - \sqrt{1+\alpha^2})$, the procedure above
yields the solutions:\\
\begin{center}
\boxed{u(x,y)=-\frac{\sqrt{1-\alpha^2}}{\alpha}x^2+xy+g\left
    (y-\frac{\sqrt{1-\alpha^2}}{\alpha}x\right)}\\
\end{center}

\begin{figure}
\centering 
\mbox{\subfigure[Example \ref{hel}, $g(s)=\tan^{-1}(s)$]{\epsfig{figure=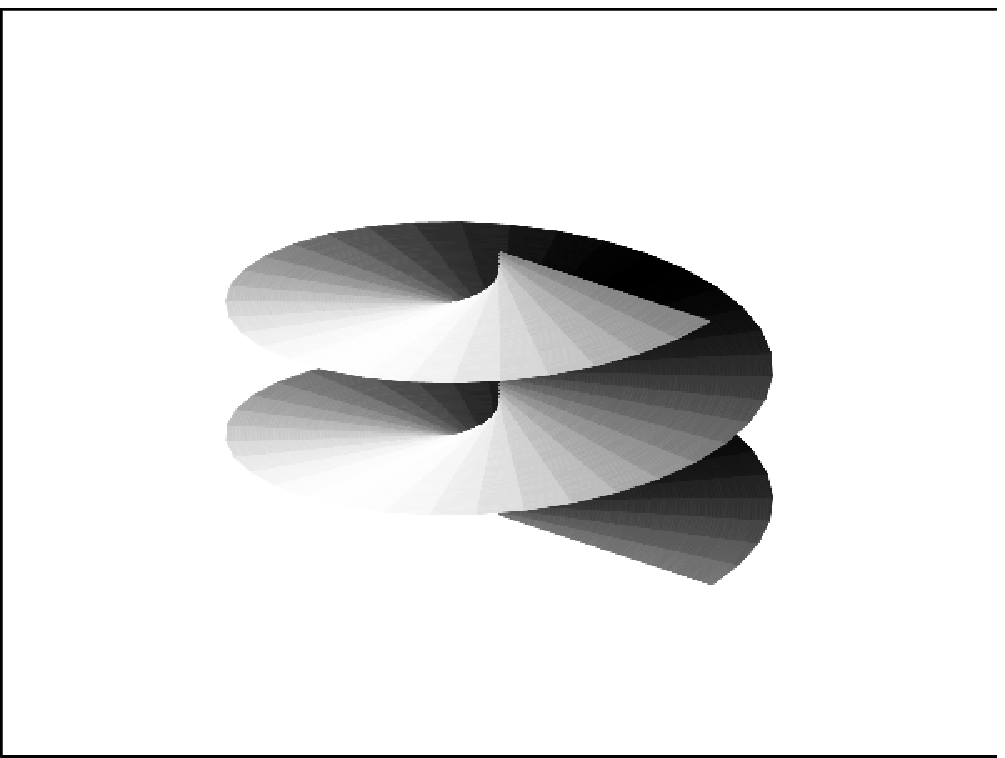,width=2.7in}}\quad
\subfigure[Example \ref{hel}, $g(s)=\frac{s}{\sqrt{1+s^2}}$]{\epsfig{figure=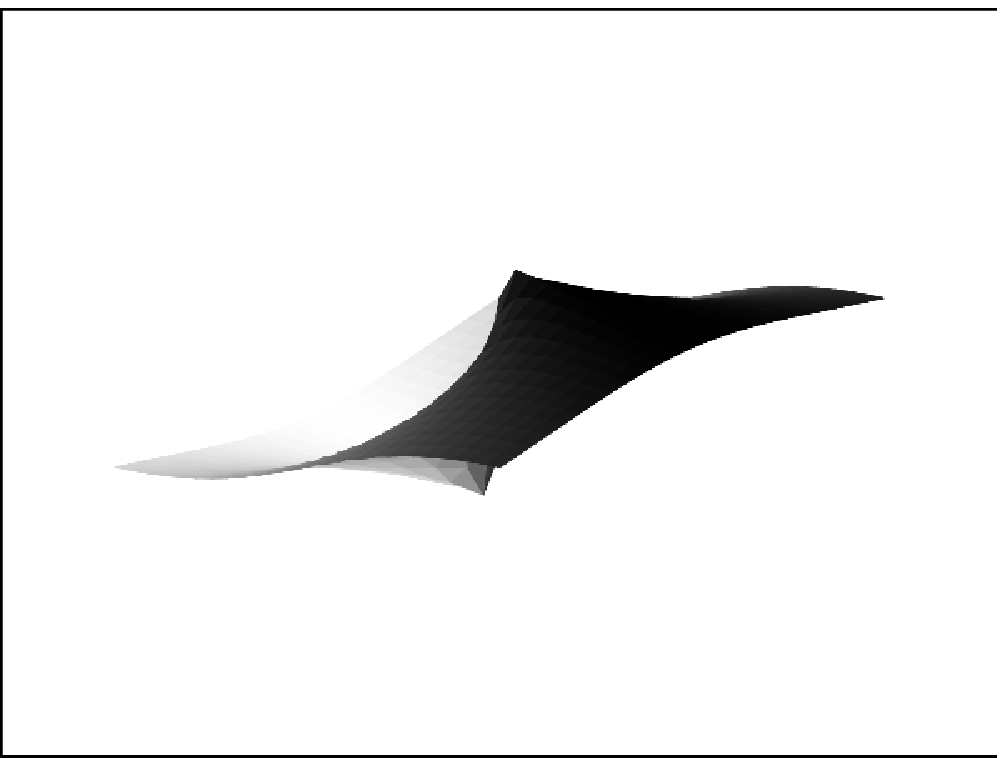,width=2.7in}}}
\mbox{
\subfigure[Example \ref{plane}, $\alpha=\beta=1, g(s)=0$]{\epsfig{figure=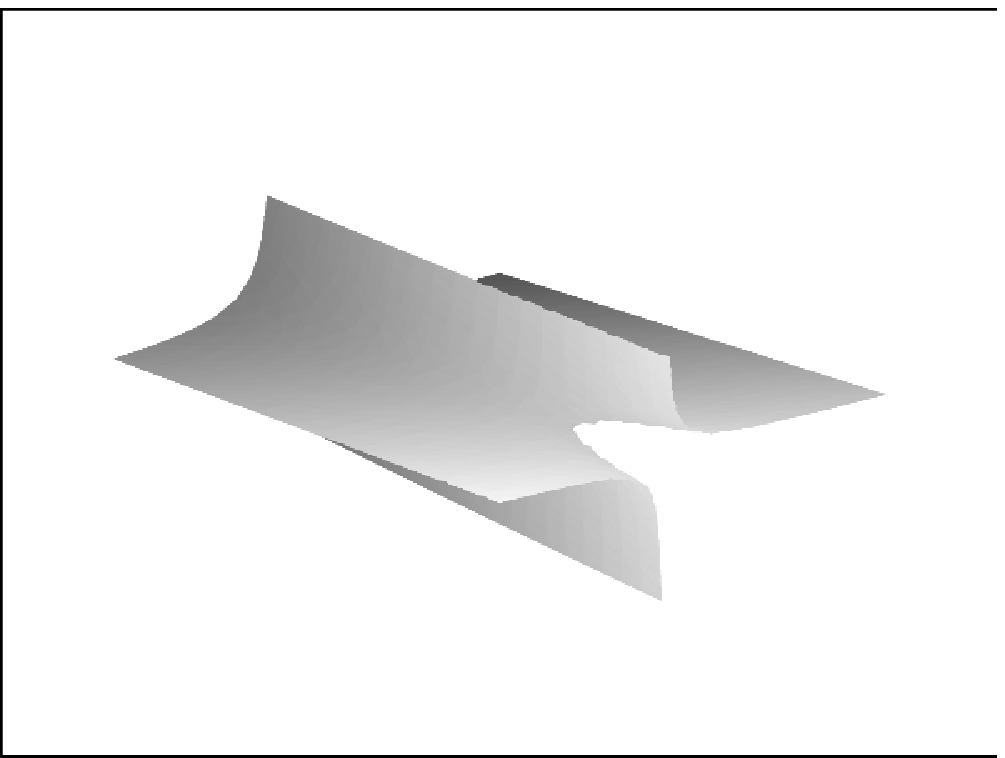,width=2.7in}}\quad
\subfigure[Example \ref{plane}, $\alpha=\beta=1, g(s)=s^2$]{\epsfig{figure=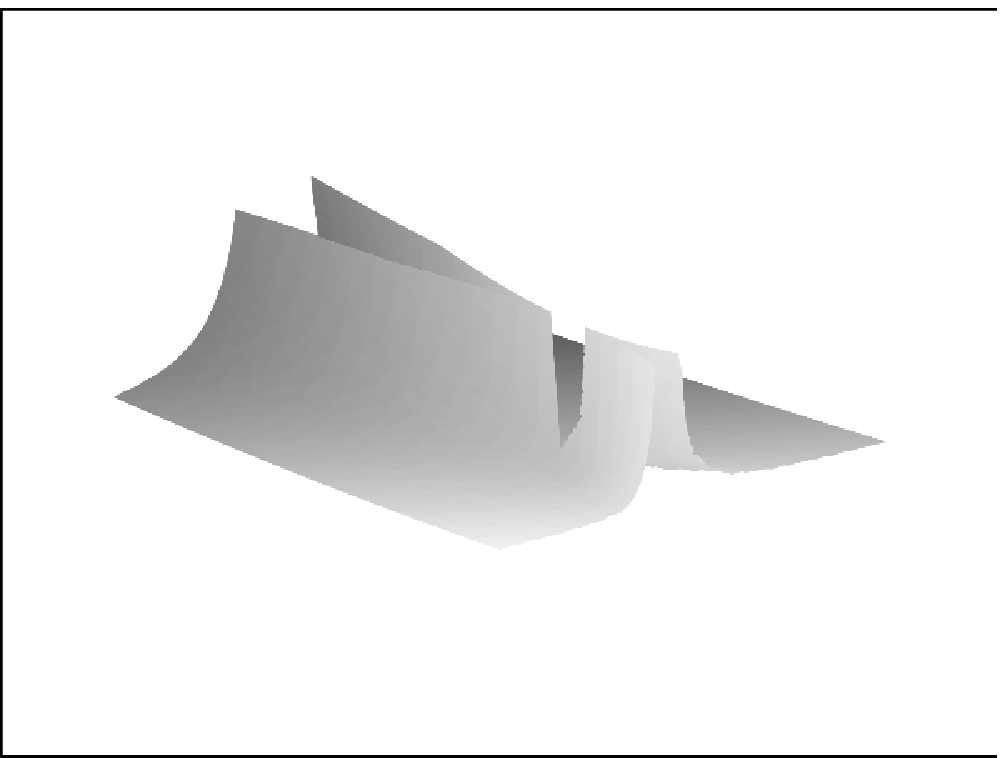,width=2.7in}}}
\caption{Minimal surfaces with Gauss map the same as some plane}\label{plfig}
\end{figure}

\subsubsection{Helicoidal examples}\label{hel}

As we saw in section \ref{inv}, both the plane $z=0$ and the helicoid
$z=\theta$ are solutions to $(\ref{MSE})$.  Just as in the previous
section, we will recover these solutions and many more from
considering the divergence free vector field associated to these
solution.  For $z=0$ or $z=\theta$, we have 
\[n_0 = \left ( \frac{-y}{\sqrt{x^2+y^2}},\frac{x}{\sqrt{x^2+y^2}}
\right)\]
Taking this vector field as $V$, $(\ref{linpde})$ becomes
\[
\frac{x}{\sqrt{x^2+y^2}}(u_x-y)+\frac{y}{\sqrt{x^2+y^2}}(u_y+x)=0
\]
or 
\begin{equation}\label{lin2}
xu_x+yu_y=0 
\end{equation}
Again, solving $(\ref{lin2})$ yields a family of solutions:\\
\begin{center}
\boxed{u(x,y)=g\left (\frac{y}{x}\right)}\\
\end{center}
where $g \in W^{1,p}(\R)$.  Note that with $g(s)=0$ and
$g(s)=\tan^{-1}(s)$, we recover $z=0$ and $z=\theta$ respectively.
Figure \ref{plfig}a,b give two examples of this type with specific choices
of $g$.  As in figure \ref{constN}

\subsubsection{Solutions associated to planes}\label{plane}
Expanding on the results of the previous section, we consider the
plane $z=a x+b y+c$ which has horizontal Gauss map
\[N_0=\left ( \frac{a -y}{\sqrt{(a-y)^2+(b+x)^2}},
  \frac{b+x}{\sqrt{(a-y)^2+(b+x)^2}} \right)\]
Again using this as $V$, $(\ref{linpde})$ becomes
\[ (b+x)(u_x-y)-(a-y)(u_y+x)=0\]
Solving this first order partial differential equation yields
solutions:\\
\begin{center}
\boxed{u(x,y)= \frac{a}{b+x} x^2+\frac{b^2}{b +x} xy
  + g\left (\frac{y-a}{b+x}\right)}\\
\end{center}
where, again, $g \in W^{1,p}(\R)$.  Figures \ref{plfig}c,d shows some
examples with different choices of $g$ and $a=b=1$.  In these figures, the shading
is not connected to the horizontal Gauss map.  Notice
that if we take $g(s)=s$, we recover the plane $z=x+y-1$.  

\subsection{Summary}
Briefly, we revisit the examples generated above, listing them
together:
\begin{Thm} \label{list} Given $a,b,c \in \R$, $\alpha \in [-1,1]$
  and $g \in W^{1,p}(\R)$, the following graphs satisfy (\ref{MSE})
  weakly:
\begin{enumerate}
\item \label{s1} \[z=a x+b y+c\]
\item \label{s2}\[z=a \theta\] (using cylindrical coordinates)
\item \label{s3}\[z= \pm \left ( \frac{\sqrt{b r^2-1}}{b} - a
      \tan^{-1} \left ( \frac{1}{\sqrt{b r^2-1}} \right ) \right )
    + c + a \theta\]  Again, these use cylindrical
    coordinates.
\item \label{s4}\[z=xy+g(y)\]
\item \label{s5}\[ z = \frac{\alpha}{\sqrt{1-\alpha^2}} x^2+xy+g\left(
    y+\frac{\alpha}{\sqrt{1-\alpha^2}}x\right)\]
\item \label{s6}\[ z = -\frac{\sqrt{1-\alpha^2}}{\alpha} x^2+xy+g\left(
    y-\frac{\sqrt{1-\alpha^2}}{\alpha}x\right)\]
\item \label{s7}\[z=g \left ( \frac{y}{x}\right)\]
\item \label{s8}\[z=\frac{a}{b+x}x^2+\frac{b^2}{b+x}xy+g \left
    ( \frac{y-a}{b+x} \right )\]
\end{enumerate}
\end{Thm}

\section{Consequences and Discussion}\label{cons}
The examples of the last sections show that there are a surprising
number of minimal surfaces in the setting of the Carnot Heisenberg
group.  In particular, the examples with prescribed Gauss map are
particularly enlightening.  These show that even when restricting the
horizontal Gauss map to give a specific vector field (up to the
ambiguity of sign), one can produce an envelope of infinitely many
solutions.  This is in stark contrast to the Euclidean and Riemannian
cases where, under suitable conditions, one expects reasonable rigidity
of minimal surfaces.  In the next two sections, we explore the
implications of the examples above in the context of two classical
problems.  The first, the so-called Bernstein problem, asks for the
classification of complete minimal graphs.  The second deals with the
uniqueness of solutions to the Dirichlet problem in this setting.  

\subsection{The Bernstein Problem}\label{Bernprob} In the theory of minimal surfaces
in Euclidean space, Berstein characterized complete minimal graphs:
\begin{Thm}[Bernstein \cite{Bernstein}] Planes are the only complete minimal graphs in $\R^3$.
\end{Thm}

One can prove this beautiful theorem by showing that the Gauss map can be
viewed as a bounded holomorphic function on $\C$ and thus must be
constant.  Therefore, the normal to the complete minimal graph is
constant, i.e. the graph is a plane.  To mark the contrast between
minimal surface theory in $\R^3$ and in $(H,\cc)$, we point out that
many of the examples listed in section \ref{egs} are complete minimal
graphs.  In particular we've seen that we can even find many complete
minimal graphs with piecewise constant Gauss map.  Theorem \ref{list}
(\ref{s4}) (\ref{s5}) and (\ref{s6}) are examples of this form.  We
emphasize that in most of these examples, one has a {\em huge} amount
of flexibility  - examples of the form (\ref{s5}) and (\ref{s6})
depend on a real number as well as an element of $W^{1,p}(\R)$.  This,
of course, shows that nothing like Bernstein's theorem can possibly be
true in this setting.  

However, if one does not allow discontinuities in the horizontal Gauss
map, the question remains:  do there exist minimal graphs with
constant horizontal Gauss map?  If so, how many are there?  Having
continuous constant Gauss map implies that the surface may not have
any characteristic points.  Since the submission of this article the
author and N. Garofalo have shown (\cite{GP02}) a version of the
Bernstein property in the Heisenberg group.

\subsection{Uniqueness for the Dirichlet problem}\label{unique}
In this section, we investigate the uniqueness of solutions to the
Dirichlet problem.  To be precise, we will investigate the following
problem:
\begin{Prob}\label{pp}  Given a set $\Omega \subset \R^2$ with a
  smooth closed boundary and a smooth closed curve $\Gamma:\partial
  \Omega \ra H$,
  does there exist $u \in W^{1,p}(\Omega)$ so that the graph
  $z=u(x,y)$ defines a surface satisfying the minimal surface equation bounded by $\Gamma$.  
\end{Prob}

In deriving the minimal surface equation $(\ref{MSE})$ in the third
section, we noted that the equation is sub-elliptic and hence, we
cannot immediately guarantee uniqueness of solutions to the Dirichlet
problem.  Indeed, this observation coupled with the wealth of examples
produced in section \ref{egs} leads one to suspect that solutions may
not be unique.  In this section we produce two distinct solutions to
the Dirichlet problem for a given curve by exploiting the connection
between surfaces satisfying (\ref{MSE}) and sets of minimal perimeter
given in theorem \ref{min2}.  We first define the minimal
surfaces themselves:

\begin{equation*}\label{g1}
u_1(x,y)=x^2+xy \tag{$S_1$}
\end{equation*}
Note that this is an example of the form produced in section
\ref{const} or theorem \ref{list} (\ref{s5}) with $\delta=\frac{\sqrt{2}}{2}$ and $f(s)=0$.

\begin{equation*}\label{g2}
u_2(x,y)=xy+1-y^2 \tag{$S_2$}
\end{equation*}
Note that this is an example of the form produced in section
\ref{const} or theorem \ref{list} (\ref{s5})  with $\delta=0$ and $f(s)=1-s^2$.

\begin{figure}[h]\label{diri}
\epsfig{file=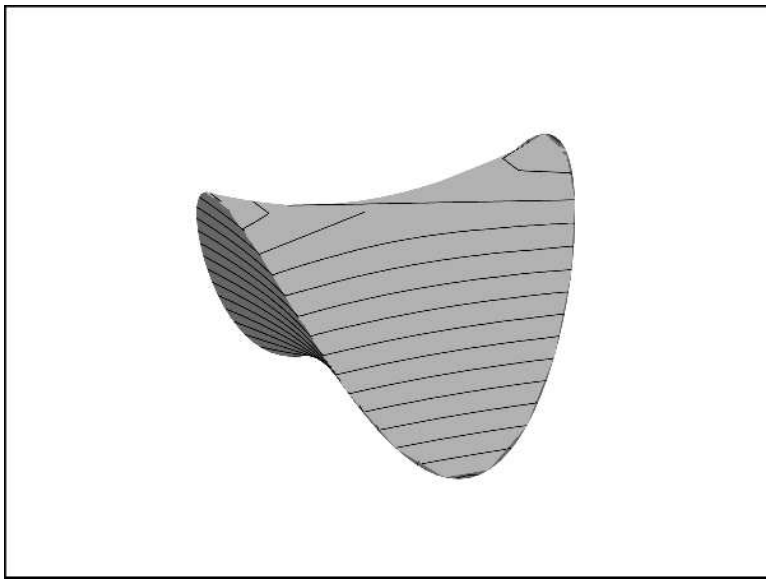}
\caption{Two minimal surfaces bounded by the same curve.}
\end{figure}

Observe
next that both surfaces contain the same smooth closed curve, $\Gamma$,
which is the graph over the unit circle in the plane given
parametrically as
\[ \Gamma=
\{(\cos(\theta),\sin(\theta),\cos^2(\theta)+\sin(\theta)\cos(\theta))
| \theta \in [0,2\pi]\}\]
Defining $\phi(\theta)=\cos^2(\theta)+\sin(\theta)\cos(\theta)$,
$\Gamma$ is the graph of $\phi$ over $S^1$.  Figure \ref{diri} shows
both curves bounded by $\Gamma$.

We summarize this observation as follows:
\begin{Thm}
The graphs given by $z=u_1(x,y)$ and $z=u_2(x,y)$ over the closed unit
ball in $\R^2$ are solutions to the Dirichlet problem:
\begin{gather*}
H_{cc}(u)=0 \\
u|_{S^1} = \phi
\end{gather*}

\end{Thm}

\noindent
{\em Remark:}  Notice that a direct calculation verifies that the area
of the two surfaces is the same:

\[ E(u_1)= \int_\Omega 8x^2 \; dx dy = 2\pi\]
and
\[ E(u_2)=\int_\Omega 4(x-y)^2 \; dx dy = 2\pi\]

\subsection{Further investigations}
The wealth of examples of minimal surfaces described above, while
answering several questions concerning the geometry of the Carnot
Heisenberg group, leads us to many unsolved problems.  In this
section we outline a few of these questions:
\renewcommand{\theenumi}{\arabic{enumi}}

\begin{enumerate}
\item In section \ref{Bernprob}, we mentioned the question of the
  existence of examples of graphs with constant Gauss map.  More generally, we
  can pose the following question:
\begin{Qu}  Can one completely classify all examples of prescribed
  piecewise constant Gauss map?
\end{Qu}

\item While the examples show that without any extra conditions, one
  cannot hope for the higher regularity of solutions, the question of
  higher regularity for solutions of the Dirichlet problem is not
  ruled out by the examples.  
\begin{Qu}
Given a closed curve, $\Gamma$, in $H$ of class $C^k$ given as a graph
over a closed curve in $\R^2$, what is the regularity of
a minimal surface spanning $\Gamma$?
\end{Qu}
\item This paper focused on nonparametric minimal surfaces, mostly
  because the calculations were not nearly as cumbersome as in the
  general case.  However, the case explored here, that of graphs over
  the $xy$-plane is set up precisely to match the nonisotropic
  character of the Carnot-Carath\'eodory metric.  It is not hard to
  write down an analogous set of equations for graphs over the
  $yz$-plane, but the resulting partial differential equation is much
  less tractable and it is not clear that one recovers $X$-minimal
  surfaces from the solutions.  
\begin{Qu}
Can one extend these techniques to find $X$-minimal surfaces which are
not necessarily graphs over domains in the $xy$-plane?
\end{Qu}

\item Understanding the case of $H_{cc}=0$ is a first step in the
  exploration of the analogue of constant mean curvature
  surfaces in the Carnot-Carath\'eodory setting.  Naively, one can
  simply find examples, using the same techniques as above, to
  solutions of $H_{cc}=const$.  For example, one can find a closed
  solution to $H_{cc}=1$, the top half of which is given by the function
\[ f(r)=-\frac{1}{2}r\sqrt{4-r^2}+2 \arctan\left
  (\frac{r}{\sqrt{4-r^2}}\right) \]
where, like the rotationally invariant solutions detailed above, we
use polar coordinates to describe the graph.  Figure \ref{egg} shows a
picture of the whole surface.  

\begin{figure}[h]
\epsfig{file=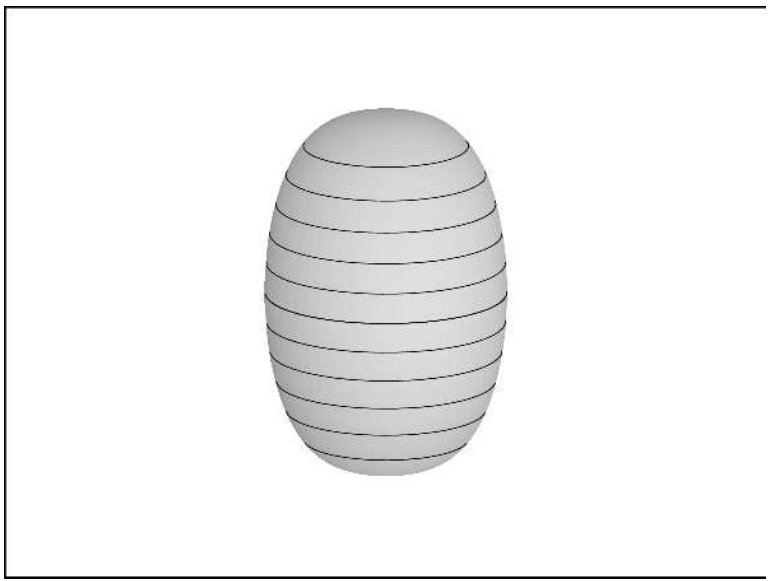}
\caption{A closed suface satisfying $H_{cc}=1$}\label{egg}
\end{figure}

The solutions of this equation do indeed minimize area among
candidates of fixed volume and so are an appropriate extension of the
notion of constant mean curvature surfaces in the Euclidean setting.

In
\cite{DGN}, Danielli, Garofalo and Nhieu address the notion of
``constant X-mean curvature surfaces'' and use them to investigate the
isoperimetric profile of $(H,\cc)$.  Specifically, they find that the
surfaces giving the best isoperimetric constant among rotationally
symmetric surfaces are indeed those satisfying the equation
$H_{cc}=const$.  This is a tantalizing result which provides a window
into the possibilities for understanding completely the connection
between the isoperimetric profile and this class of surfaces.  In
particular, it gives hope that, as in many Riemannian settings, the
constant mean curvature surfaces reflect the isoperimetric profile
exactly.

\end{enumerate}

%\bibliographystyle{alpha}
%\bibliography{../biblio}

\end{document}